\documentclass[11pt,a4paper,reqno]{amsart}

\usepackage{amssymb,enumerate}
\usepackage{mathtools}
\usepackage{microtype}
\usepackage[myheadings]{fullpage}
\usepackage{tikz}
\usepackage{graphicx}
\usepackage{subcaption}
\usepackage{hyperref}
\hypersetup{colorlinks=true,linkcolor=blue,citecolor=blue}

\DeclareMathOperator{\Stab}{Stab}
\DeclareMathOperator{\Aut}{Aut}
\DeclareMathOperator{\Out}{Out}
\newtheorem{thm*}{Theorem}
\newtheorem{cor*}{Corollary}
\newtheorem{theorem}{Theorem}[section]
\newtheorem{cor}[theorem]{Corollary}
\newtheorem{prop}[theorem]{Proposition}
\newtheorem{lemma}[theorem]{Lemma}
\theoremstyle{definition}
\newtheorem{definition}[theorem]{Definition}
\newtheorem{rem}[theorem]{Remark}

\begin{document}

\title[Liftable outer automorphisms of cyclic covers]{Generating the liftable outer automorphism groups of regular cyclic covers}
\author{Rohan Suresh Mahure}
\address{(R. S. Mahure) Department of Mathematics\\
Indian Institute of Science Education and Research Bhopal\\
Bhopal Bypass Road, Bhauri \\
Bhopal 462 066, Madhya Pradesh\\
India}
\email{rohansmahure@gmail.com, mahure22@iiserb.ac.in}
\urladdr{https://sites.google.com/iiserb.ac.in/rohanmahure/home}

\author{Priyavrat Pandey}
\address{(P. Pandey) Department of Mathematics\\
Indian Institute of Science Education and Research Bhopal\\
Bhopal Bypass Road, Bhauri \\
Bhopal 462 066, Madhya Pradesh\\
India}
\email{priyavrat23@iiserb.ac.in}

\author{Kashyap Rajeevsarathy}
\address{(K. Rajeevsarathy) Department of Mathematics\\
Indian Institute of Science Education and Research Bhopal\\
Bhopal Bypass Road, Bhauri \\
Bhopal 462 066, Madhya Pradesh\\
India}
\email{kashyap@iiserb.ac.in}
\urladdr{https://home.iiserb.ac.in/~kashyap/}

\subjclass[2020]{Primary 57K20; Secondary 57M60}

\keywords{free groups, outer automorphisms, cyclic covers}

\begin{abstract}
Let $R_n = \bigvee_{i=1}^{n} S^1$ be a finite connected one-dimensional CW complex with $n \geq 2$, and let $\mathrm{Out}(F_n)$ be the outer automorphism group of the free
group $F_n$, where $\pi_1(R_n)\cong F_n$. Let $p_k:\tilde{R}_n\rightarrow R_n$ be the standard $k$-sheeted regular cyclic cover of $R_n$ induced by a $\mathbb{Z}_k$-action on $\tilde{R}_n$ (which is homotopy equivalent to $R_{n_k}$, where $n_k=k(n-1)+1$) generated by a free rotation through $2\pi/k$ of $\tilde{R}_n$. Let $\mathrm{Out}_{F_{n_k}}(F_n)$ be the liftable outer automorphism group associated with a standard $k$-sheeted normal cover $p_k$. For $n\geq3$, we show that $\{\mathrm{Out}_{F_{n_k}}(F_n)\}_{k\geq2}$ forms an infinite family of self-normalizing subgroups of $\mathrm{Out}(F_n)$. For $n\geq2$, these groups are maximal if and only if $k$ is prime. Furthermore, we derive explicit finite generating sets for $\mathrm{Out}_{F_{n_k}}(F_n)$ for $k\geq 2$ and $n\geq 3$.  Let $\mathcal{C}$ be the infinite chain of roses, and let $q_n:\mathcal{C}\longrightarrow R_n$ be the standard infinite-sheeted cover induced by $\langle \tau\rangle$, where $\tau$ is the unit translation on $\mathcal{C}$. As a final application, we obtain a finite generating set for $\mathrm{Out}_{F_\infty}(F_n)$ when $n\geq 3$.
\end{abstract}

\maketitle

\section{Introduction}\label{sec:intro}
For $n\ge 2$, let $F_n=\langle x_1,\ldots,x_n\rangle$ be the free group of rank $n$. Let $R_n=\bigvee_{i=1}^{n}S^1$ denote the rose with $n$ petals, so that $\pi_1(R_n)\cong F_n$, and let $\Out(F_n)$ be the outer automorphism group of $F_n$. Together with the mapping class group of a surface, $\Out(F_n)$ is one of the central objects of geometric group theory, and much of the work on it is guided by the analogy between the two groups. Under this analogy, Outer space plays the role of Teichmüller space; the homology representation $\Psi\colon \Out(F_n)\to\mathrm{GL}(n,\mathbb{Z})$, induced by abelianizing $F_n$, plays the role of the symplectic representation; and the kernel of $\Psi$, known as the Torelli subgroup, plays the role of the Torelli subgroup of mapping class groups. Within this framework, a natural question concerns finite-index subgroups of $\Out(F_n)$, their algebraic structure, and how they are embedded in the ambient group. A special class of finite-index subgroups arises as the \emph{liftable} subgroups of finite-sheeted covers. The study of liftable self-homotopy equivalences goes back to Sasao and Cha~\cite{Sasao}, and was refined for $\Out(F_n)$ by Bridson and Vogtmann~\cite{Bridson}. More recently, the groups have been analyzed for finite-index subgroups $H\le F_n$ and finitely generated subgroups $\Gamma\le\Out(F_n)$ in~\cite{spence}. The relationship between mapping class groups and automorphism groups of free groups has been an important theme in geometric group theory (see~\cite{brid_survey}). In particular, the fundamental groups of connected finite-type surfaces with nonempty boundary or at least one puncture are free groups of finite rank, providing a natural connection between mapping class groups and $\Out(F_n)$. In the setting of surfaces, the liftable mapping class groups of regular cyclic covers were studied by Agarwal et al.~\cite{Agarwal_lift}, who obtained a symplectic characterization of their elements, and the generating sets and further structural properties of these groups were then investigated by Dey et al.~\cite{Rajeev}. Motivated by these works, we ask the analogous questions for $\Out(F_n)$, adapting the group-theoretic methods of~\cite{Agarwal_lift,Rajeev} where they carry over and developing new arguments where they do not.

\begin{figure}[htbp]
    \centering
   \begin{subfigure}[b]{0.48\textwidth}
        \centering
        \scalebox{0.78}{%
\begin{tikzpicture}
  \node[inner sep=0] (fig) at (0,0) {\includegraphics[width=5cm]{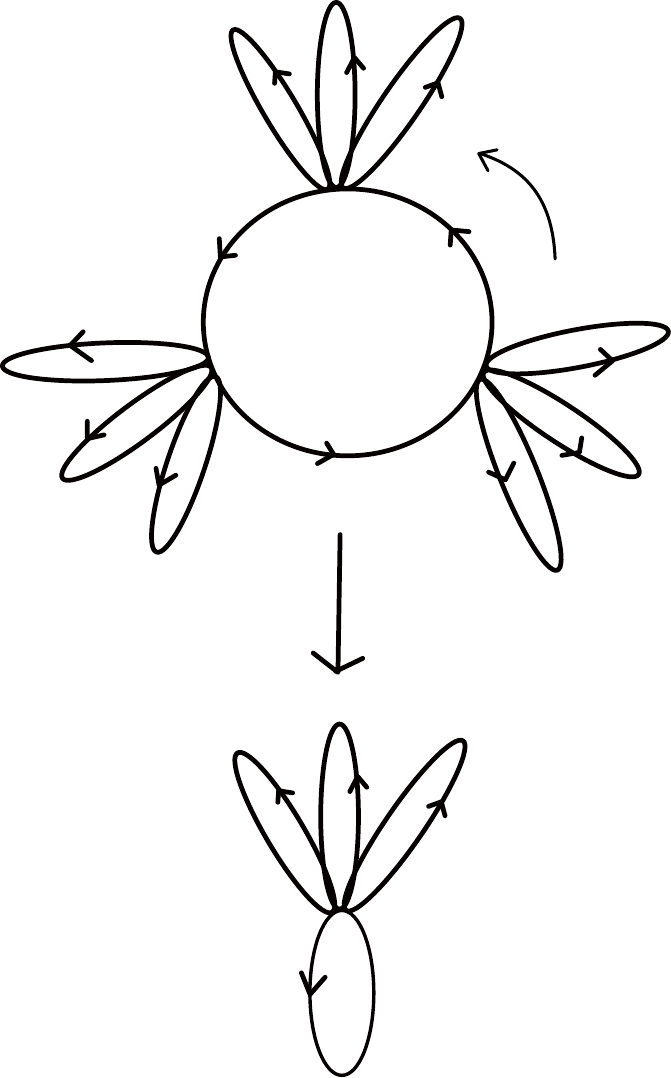}};
  \node at (.9,1.35)  {$v_1$};
  \node at (.15,2.4) {$v_2$};
  \node at (-.65,1.35)  {$v_3$};
  \node at (-1,3.6) {$a_4^2$};
  \node at (0.9,4.2)  {$a_2^2$};
  \node at (-0.2,4.3)  {$a_3^2$};
  \node at (2.7,.7) {$a_3^1$};
  \node at (-1.6,-0.4) {$a_4^3$};
  \node at (-2.4,.4) {$a_3^3$};
  \node at (2.8,1.6)  {$a_4^1$};
  \node at (2.0,-0.3) {$a_2^1$};
  \node at (-2.8,1.4) {$a_2^3$};
  \node at (0,.27) {$ a_1$};
  \node at (2,3)  {$2\pi/3$};
  \node at (0.35,-0.5) {$p_3$};
  \node at (-0.5,-3.5) {$\bar a_1$};
  \node at (1.1,-1.9) {$\bar a_2$};
  \node at (-.25,-1.4)  {$\bar a_3$};
  \node at (-.93,-1.9)  {$\bar a_4$};
\end{tikzpicture}%
}
\caption{A regular cyclic cover $p_3:\tilde{R}_4\to R_4$ induced by a $\mathbb Z_3$-action on $\tilde{R}_4$.}
\label{fig:z3-cover}
    \end{subfigure}
    \hfill
   \begin{subfigure}[b]{0.48\textwidth}
        \centering
        \scalebox{0.78}{%
\begin{tikzpicture}
  \node[inner sep=0] (fig) at (0,0) {\includegraphics[width=9.8cm]{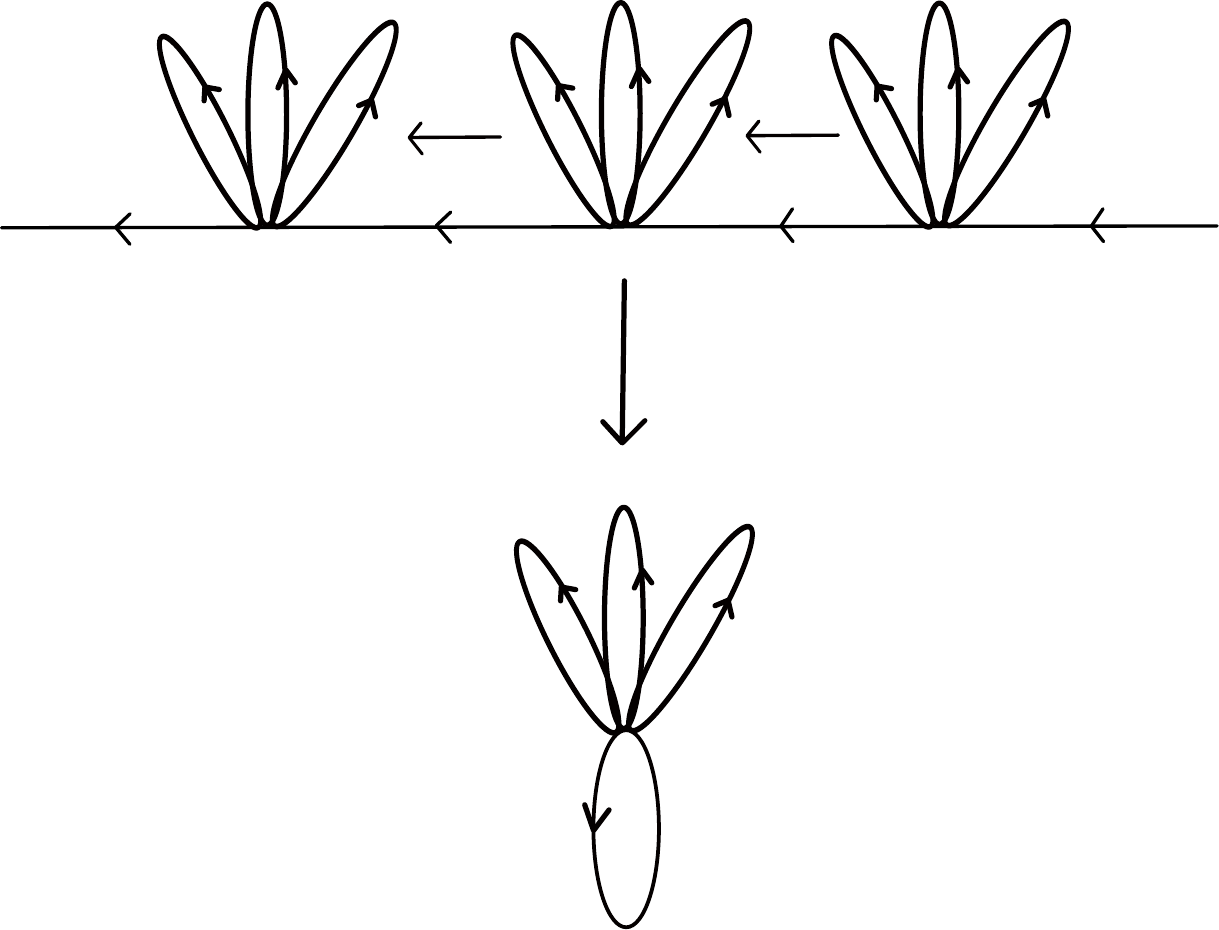}};
  \node at (-1.2,2.9) {$ \tau$};
\node at (.9,2.2) {$ a_2^{0}$};
  \node at (.3,3.9)  {$ a_3^{0}$};
\node at (-.7,2.2)  {$ a_4^{0}$};
\node at (-1.8,2.2) {$ a_2^{-1}$};
  \node at (-2.4,3.9)  {$ a_3^{-1}$};
\node at (-3.5,2.2)  {$ a_4^{-1}$};
 \node at (3.4,2.2) {$ a_2^{1}$};
  \node at (2.9,3.9)  {$ a_3^{1}$};
\node at (1.8,2.2)  {$ a_4^{1}$};
 \node at (2.6,1.65) {$ v_1$};
  \node at (.2,1.65)  {$ v_0$};
\node at (-2.5,1.65)  {$ v_{-1}$};
  \node at (-.43,-2.8) {$\bar a_1$};
  \node at (1.26,-.35) {$\bar a_2$};
  \node at (-0.26,-.35)  {$\bar a_3$};
\node at (-.7,-1.4)  {$\bar a_4$};
\node at (-1.5,1.6)  {$ e_{1,0}$};
\node at (1.5,1.6)  {$ e_{1,1}$};
 \node at (.4,.9)  {$ q_4$};
\end{tikzpicture}%
}
\caption{A regular cyclic cover $q_4:\mathcal{C}\to R_4$ induced by a $\mathbb{Z}$-action on $\mathcal{C}$.}
\label{fig:inf-cover}
    \end{subfigure}
    \caption{Regular cyclic covers of $R_4$ by $\tilde{R}_4$ and $\mathcal{C}$.}
    \label{fig:cyclic-covers}
\end{figure}

For $k\ge2$, let $p_k:\tilde{R}_n\to R_n$ be the standard regular $k$-sheeted cover of $R_n$ induced by a $\mathbb Z_k$-action on $\tilde{R}_n$ generated by a free rotation through the angle $2\pi/k$, as shown in Figure~\ref{fig:z3-cover} (for $k=3$). The vertex set of $\tilde{R}_n$ is $\{v_1,\dots,v_k\}$, with one vertex for each sheet, and the edges $v_j\to v_{j+1\,(\mathrm{mod}\,k)}$ form the central circle $C_k$, with $n-1$ petals attached at each vertex. In particular, $\tilde{R}_n$ is homotopy equivalent to $R_{n_k}$, where $n_k=k(n-1)+1$. Let $\Out_{F_{n_k}}(F_n)$ denote the \emph{liftable} subgroup of $\Out(F_n)$ associated with $p_k$, that is, the group of outer classes that are represented by self-homotopy equivalences of $R_n$ lifting to $\tilde{R}_n$.  Since there are only finitely many epimorphisms $H_1(R_n,\mathbb Z)\to\mathbb Z_k$, we have $[\Out(F_n):\Out_{F_{n_k}}(F_n)]<\infty$, and consequently $\Out_{F_{n_k}}(F_n)$ is finitely generated. In this paper, our objective is to analyze this subgroup by characterizing it homologically, identifying an explicit finite generating set, and examining its properties.
Consider the standard generators $\{a_i^j: 2\le i\le n,\ 1\le j\le k\}\cup\{a_1\}$ of $H_1(\tilde{R}_n,\mathbb Z)$ and $\{\bar a_i:1\le i\le n\}$ of $H_1(R_n,\mathbb Z)\cong\mathbb Z^n$, where $a_1$ represents the central circle, and let $e_1,\dots,e_n$ denote the standard unit vectors representing $\bar a_1,\dots,\bar a_n$. By identifying $\pi_1(\tilde{R}_n)$ with its image under $(p_k)_*$, the induced map on homology is
\begin{equation}\label{eq:finite-cover-homology}
(p_k)_\#(a_i^j)=\bar{a}_i
\quad (2\leq i\leq n,\ 1\leq j\leq k),
\qquad
(p_k)_\#(a_1)=k\bar{a}_1.
\end{equation}

By the universal coefficient theorem~\cite{Hatcher2002},
\[
H^1(R_n,\mathbb Z_k)\cong\operatorname{Hom}(H_1(R_n,\mathbb Z),\mathbb Z_k)
\cong(\mathbb Z_k^n)^*,
\]
with the standard dual basis $\{e_1^*,e_2^*,\ldots,e_n^*\}$,
where $e_i^*(e_j)=\delta_{ij}$ for all $1\leq i,j\leq n$.
Let $ \Psi : \Out(F_n) \longrightarrow \mathrm{GL}(n,\mathbb{Z})$ be a linear representation induced by abelianizing the free group $F_n$. Nielsen showed that the map $\Psi$ is surjective~\cite{Nielsen1924}. Let $\mathcal J_{n,k}=\rho_k(\operatorname{GL}(n,\mathbb Z))$, where $\rho_k$ denotes reduction modulo $k$. Consider the surjective homomorphism
$\Psi_k=\rho_k\circ\Psi:\Out(F_n)\twoheadrightarrow
\mathcal J_{n,k}$. The group $\Out(F_n)$ acts transitively on the primitive cohomology classes in both $H^1(R_n,\mathbb Z)$ and $H^1(R_n,\mathbb Z_k)$. This transitivity allows us to express the index of $\Out_{F_{n_k}}(F_n)$ as an orbit count.

In Section~\ref{sec:homological_criteria}, we apply basic covering space theory to obtain a complete characterization of $\Out_{F_{n_k}}(F_n)$ in terms of its image under both $\Psi$ and $\Psi_k$ as stated in the following theorem.
\begin{thm*}\label{thm:intro-liftability}
Given $\phi \in \Out(F_n)$, the following statements are equivalent.
\begin{enumerate}[(i)]
\item $\phi \in \Out_{F_{n_k}}(F_n)$.
\item If $\Psi(\phi) = (c_{ij})_{n \times n}$, then $k\mid c_{1j}$ for $2\leq j\leq n$ and $\gcd(c_{11},k)=1$.
\item  If $\Psi_k(\phi) = (e_{ij})_{n\times n}$, then $e_{11}\in\mathbb Z_k^\times$ and $e_{1j}=0$ for $2\leq j\leq n$.
\item $\phi\in\Stab_{\Out(F_n)}(\{\ell e_1^*:\ell\in\mathbb Z_k^\times\})$.
\end{enumerate}
\end{thm*}
\noindent The following corollaries describe intersections and indices of liftable subgroups and follow from Theorem~\ref{thm:intro-liftability}.
\begin{cor*}\label{cor:intro-intersections}
For $n \geq 2$ and $k, \ell \geq 2$, we have
\begin{enumerate}[(i)]
  \item $\mathcal{T}_{n,k} \triangleleft \Out_{F_{n_k}}(F_n)$.
  \item $\Out_{F_{n_k}}(F_n)\cap\Out_{F_{n_\ell}}(F_n)=\Out_{F_{n_d}}(F_n),$ where $d=\operatorname{lcm}(k,\ell)$.
\end{enumerate}
\end{cor*}
\begin{cor*}
Let $n\geq 2$ and $k\geq 2$. Let $\mathcal{I}_{n,k}^{\mathbb Z}=\Psi(\Out_{F_{n_k}}(F_n))$ and
$\mathcal I_{n,k}=\Psi_k(\Out_{F_{n_k}}(F_n))$. Let $\mathcal J_{n,k}=\rho_k(\mathrm{GL}(n,\mathbb Z))$, where $\rho_k$ denotes reduction modulo $k$. Then the following statements hold.
\begin{enumerate}[(i)]
\item
$[\Out(F_n):\Out_{F_{n_k}}(F_n)]
=[\mathrm{GL}(n,\mathbb{Z}):\mathcal{I}_{n,k}^{\mathbb Z}]
=[\mathcal{J}_{n,k}:\mathcal{I}_{n,k}]$.
\item We have
\[
[\Out(F_n):\Out_{F_{n_k}}(F_n)]
=\frac{|\{\text{primitive elements of }(\mathbb Z_k^n)^*\}|}{E(k)}.
\]
In particular, if $k$ is prime, then
\[
[\Out(F_n):\Out_{F_{n_k}}(F_n)]=\frac{k^n-1}{k-1}.
\]
\end{enumerate}
Here $E(k)$ is Euler's totient function.
\end{cor*}
Let $\mathcal{T}_n$ denote the Torelli subgroup, defined as the kernel of $\Psi$. For each $k\geq2$, let $\mathcal{T}_{n,k}$ be defined as the kernel of $\Psi_k$. The subgroup $\Out_{F_{n_k}}(F_n,e_1^*)$ is defined as the stabilizer of $e_1^*$ in $\Out(F_n)$, that is, $\Out_{F_{n_k}}(F_n,e_1^*)=\Stab_{\Out(F_n)}(e_1^*)$. In Section~\ref{sec:criter_stab}, we give homological criteria for membership in $\Out_{F_{n_k}}(F_n,e_1^*)$. For distinct indices $1\leq i,j\leq n$, let $\mu_{ij}$ denote the Nielsen transformation of $F_n$ defined by $\mu_{ij}(x_i)=x_i x_j$ and $\mu_{ij}(x_\ell)=x_\ell$ for all $\ell\neq i$.  In Section~\ref{sec:normal_series}, we derive a normal series for $\Out_{F_{n_k}}(F_n)$, as an application of Theorem~\ref{thm:homo_criteria}.
\begin{thm*}\label{thm:intro-normal-series}
Let $n\geq 2$ and $k\geq 2$. Then $\Out_{F_{n_k}}(F_n)$ admits a normal series
\[
1\triangleleft\mathcal{T}_{n,k}
\triangleleft\Out_{F_{n_k}}(F_n,e_1^*)
\triangleleft\Out_{F_{n_k}}(F_n).
\]
Furthermore,
\[
\Out_{F_{n_k}}(F_n)/
\Out_{F_{n_k}}(F_n,e_1^*)
\cong\mathbb{Z}_k^\times,
\]
and
\[
[\Out(F_n):\Out_{F_{n_k}}(F_n,e_1^*)]
=
\left|\left\{\text{primitive elements of }(\mathbb{Z}_k^n)^*\right\}\right|.
\]
Choose integer representatives of $\ell,\bar\ell\in\mathbb Z_k^\times$ with $\ell\bar\ell\equiv1\pmod{k}$. Moreover, a set of coset representatives for
$\Out_{F_{n_k}}(F_n)/
\Out_{F_{n_k}}(F_n,e_1^*)$ is given by
\[
\mathcal{S}_k''=
\left\{
w_\ell=
\mu_{12}^{1-\bar{\ell}}\mu_{21}^{-1}
\mu_{12}^{1-\ell}\mu_{21}^{\bar{\ell}}
\mid \ell\in\mathbb{Z}_k^\times
\right\}.
\]
\end{thm*}
For $n\geq2$, let $\iota\in\Out(F_n)$ be the class of the basis-inverting automorphism $x_i\mapsto x_i^{-1}$ for every $i$. Thus $\iota^2=1$ and $\Psi(\iota)=-I_n$. We obtain the following extension of the normal series in Theorem~\ref{thm:normal_series_app}.
\begin{cor*}\label{cor:norm_invo_intro}
For $n \geq 2$ and $k \geq 3$, the group $\Out_{F_{n_k}}(F_n)$ admits the normal series
\[
1
\triangleleft \mathcal{T}_{n,k}
\triangleleft \Out_{F_{n_k}}(F_n,e_1^*)
\triangleleft \langle\!\langle \Out_{F_{n_k}}(F_n,e_1^*),\iota
\rangle\!\rangle
\triangleleft \Out_{F_{n_k}}(F_n),
\]
where
$\langle\!\langle \Out_{F_{n_k}}(F_n,e_1^*),\iota
\rangle\!\rangle$ denotes the normal closure of
$\langle \Out_{F_{n_k}}(F_n,e_1^*),\iota\rangle$ in
$\Out_{F_{n_k}}(F_n)$,
\[
\langle\!\langle \Out_{F_{n_k}}(F_n,e_1^*),\iota
\rangle\!\rangle /
\Out_{F_{n_k}}(F_n,e_1^*)
\cong \mathbb{Z}_2,
\]
and
\[
\left[
\Out_{F_{n_k}}(F_n):
\langle\!\langle \Out_{F_{n_k}}(F_n,e_1^*),\iota
\rangle\!\rangle
\right]
= \frac{E(k)}{2}.
\]
\end{cor*}
In Sections~\ref{sec:maxima} and~\ref{sec:self_norm}, we present two structural results that describe the position of the liftable subgroup within $\Out(F_n)$. When $k$ is composite, $\Out_{F_{n_k}}(F_n)$ is properly contained in $\Out_{F_{n_p}}(F_n)$ for every prime divisor $p$ of $k$ (this follows from Theorem~\ref{thm:intro-liftability}), so maximality cannot hold for composite $k$. Thus, we obtain the following results.
\begin{thm*}\label{thm:intro-structure}
For $n \ge 3$ and $k \ge 2$,
\begin{enumerate}
  \item $\Out_{F_{n_k}}(F_n)$ is self-normalizing in $\Out(F_n)$; and
  \item when $k$ is prime, the liftable outer automorphism group $\Out_{F_{n_k}}(F_n)$ is a maximal subgroup of $\Out(F_n)$.
\end{enumerate}
\end{thm*}
We briefly explain how our proof of Theorem~\ref{thm:maximal} differs from the proof of the corresponding result in~\cite{Rajeev}. In~\cite{Rajeev}, maximality of the liftable mapping class group was established by studying its symplectic image in $\mathrm{Sp}(2g,\mathbb{Z})$ using coset representatives and relations among Dehn twists. Since these relations are not available in our setting, we instead use the following observation. By Theorem~\ref{thm:intro-liftability}, for $k$ prime, $\Out_{F_{n_k}}(F_n)$ is exactly the stabilizer of the line $\langle e_1^*\rangle$, which represents a point in the projective space $\mathbb P^{n-1}(\mathbb Z_k)$, under the action of $\Out(F_n)$ on $\mathbb P^{n-1}(\mathbb Z_k)$. This action is transitive. Since the stabilizer of a point in a transitive action is maximal if and only if the action is primitive, Lemma~\ref{lem:max_point_stab} reduces to proving that the action of $\Out(F_n)$ on $\mathbb P^{n-1}(\mathbb Z_k)$ is primitive.

In Section~\ref{sec:gen_set}, we determine a generating set for the stabilizer subgroup $\Out_{F_{n_k}}(F_n,e_1^*)$ appearing at the middle stage of this normal series. Before stating this, recall that the principal congruence subgroup of level $k$ in $\mathrm{SL}(2, \mathbb{Z})$, denoted by $\Gamma(k)$, is defined by
\[\Gamma(k)=\left\{
\begin{bmatrix}
a & b\\
c & d
\end{bmatrix}\in \mathrm{SL}(2,\mathbb{Z}):a,d\equiv 1\pmod{k}, \text{ and }\;b,c\equiv 0\pmod{k}\right\}.\]
Choose a finite generating set $\mathcal S_k$ for $\Gamma(k)$ (see~\cite[Lemma~12]{DiamantisKnoppMasonOSullivan2006} and~\cite[Proposition~1.17]{Stein2007}). For each $A\in\mathcal S_k$, choose a word in the matrices $\Psi(\mu_{12})^{\pm1}$ and $\Psi(\mu_{21})^{\pm1}$ representing $A$. Replace these matrices by the corresponding Nielsen automorphisms of $F_n$, fixing $x_3,\ldots,x_n$, and let $\widetilde A$ be the resulting outer class. Define $\widetilde{\mathcal S}_k=\{\widetilde A:A\in\mathcal S_k\}$. Thus $\Psi(\widetilde A)=\operatorname{diag}(A,I_{n-2})$. The choices are made at the level of Nielsen words in $\Aut(F_n)$; no extension of an arbitrary class in $\Out(F_2)$ is assumed. For distinct indices $1\leq i,j\leq n$, let $N_{ij}\in\Aut(F_n)$ be defined by $N_{ij}(x_i)=x_jx_ix_j^{-1}$ and $N_{ij}(x_m)=x_m$ for $m\neq i$. Similarly, for $2\leq j<\ell\leq n$, let $N_{1j\ell}\in\Aut(F_n)$ be the automorphism defined by $N_{1j\ell}(x_1)=x_1[x_j,x_\ell]$ and $N_{1j\ell}(x_m)=x_m$ for all $m\neq 1$. Define $\mathcal{N}_n=\{N_{1j},N_{j1}\mid 2\le j\le n\}\cup\{N_{1j\ell}\mid 2\le j<\ell\le n\}\cup\{N_{i\ell}\mid 2\le i,\ell\le n,\ i\ne\ell\}$.
Our main result, which provides generating sets for $\Out_{F_{n_k}}(F_n,e_1^*)$ and $\Out_{F_{n_k}}(F_n)$, is stated in the following theorem.
\begin{thm*}\label{thm:intro-generating-sets}
Let $k \ge 2$ and $n \ge 3$. Then
\[
\Out_{F_{n_k}}(F_n,e_1^*)
=
\left\langle \widetilde{\mathcal{S}}_k\cup\mathcal{N}_n\cup\{\iota_2,\mu_{23}\}
\cup\{\mu_{1i}\mid 2\leq i\leq n\}
\cup\{\sigma_{ij}\mid 2\leq i<j\leq n\}
\right\rangle.
\]
Consequently, we have
\[
\Out_{F_{n_k}}(F_n)
=
\left\langle \mathcal S_k''\cup \widetilde{\mathcal{S}}_k\cup\mathcal{N}_n\cup\{\iota_2,\mu_{23}\}
\cup\{\mu_{1i}\mid 2\leq i\leq n\}
\cup\{\sigma_{ij}\mid 2\leq i<j\leq n\}\right\rangle.
\]
\end{thm*}
\noindent In Section~\ref{sec:gen_set}, we first describe the image of $\Out_{F_{n_k}}(F_n,e_1^*)$ under $\Psi$ in terms of a generating set in Theorem~\ref{thm:gen_img_stab_e1} and describe the kernel of the restriction of $\Psi$ in Theorem~\ref{thm:kernel_Psi}. Together, these results establish the generating statement above. Since $\Out_{F_{n_k}}(F_n,e_1^*)$ is a normal subgroup of $\Out_{F_{n_k}}(F_n)$ in the series above, we can use this generating set to obtain a generating set for the full liftable group.

The fundamental groups of infinite-type surfaces of infinite genus are free groups of countably infinite rank, providing a natural connection between big mapping class groups and $\Out(F_\infty)$. Let $\mathcal C$ be the infinite chain of roses with vertices $v_j$ for $j\in\mathbb Z$, edges $e_{1,j}$ connecting $v_j$ to $v_{j-1}$, and loops $a_i^j$ for $2\leq i\leq n$ attached to each $v_j$. Motivated by the study of big mapping class groups, we consider the group $\Out(F_\infty)$ and its action on $R_\infty$. This provides a natural setting for investigating analogues of liftable mapping class groups in the context of outer automorphism groups of free groups. In Section~\ref{sec:Inf_cover}, we consider the standard infinite-sheeted cover $q_n:\mathcal{C}\longrightarrow R_n$ induced by $\langle \tau\rangle$, where $\tau$ is the unit translation on $\mathcal{C}$ (see Figure~\ref{fig:inf-cover}). Let $\Out_{F_\infty}(F_n)$ denote the liftable outer automorphism group associated with $q_n$. After establishing a short exact sequence analogous to the Birman--Hilden sequence for the cover $q_n$, we obtain the following results, which provide a generating set for $\Out_{F_\infty}(F_n)$.

\begin{thm*}\label{thm:uout_F_inf}
For $n\geq 2$, let $\operatorname{UOut}(F_n):=\bigcap_{k\geq 2}\Out_{F_{n_k}}(F_n)$. Let $\phi\in\Out(F_n)$ and $\Psi(\phi)=(c_{ij})_{n\times n}$. Then the following statements hold.
\begin{enumerate}[(i)]
\item $\phi\in\operatorname{UOut}(F_n)$ if and only if $c_{1i}=0$ for $2\leq i\leq n$ and $c_{11}=\pm1$.
\item For $n\geq 3$,
$\operatorname{UOut}(F_n)=\left\langle\mathcal{N}_n\cup\{\iota_2,\mu_{23},\iota\}\cup\{\mu_{1i}\mid 2\leq i\leq n\}\cup\{\sigma_{ij}\mid 2\leq i<j\leq n\}\right\rangle.$
\item Let $\Out(F_n,e_1^*)$ denote the stabilizer of $e_1^*\in H^1(R_n,\mathbb Z)$. There is a normal series
$1\triangleleft\mathcal{T}_n\triangleleft\Out(F_n,e_1^*)\triangleleft\operatorname{UOut}(F_n)$
in $\Out(F_n)$, where the two distinct cosets of $\operatorname{UOut}(F_n)/\Out(F_n,e_1^*)$ are represented by the elements of $\langle\iota\rangle$.
\end{enumerate}
\end{thm*}
\noindent This theorem provides a complete characterization of the liftable outer automorphism group associated with the infinite-sheeted cover $q_n$ and gives an explicit generating set for it. Moreover, when the indices $k\geq 2$ are ordered by divisibility, the liftable subgroups form an inverse system under the inclusions $\Out_{F_{n_\ell}}(F_n)\hookrightarrow\Out_{F_{n_k}}(F_n)$ for $k\mid\ell$, and $\Out_{F_\infty}(F_n)$ is the inverse limit of this system, as explained at the end of Section~\ref{sec:Inf_cover}.

The paper is organized as follows. In Section~\ref{sec:preliminaries}, we introduce the necessary background and notation. In Section~\ref{sec:homological_criteria}, we derive homological criteria for determining when an outer automorphism is liftable with respect to a regular cyclic cover. In Section~\ref{sec:criter_stab}, we derive homological criteria for determining when an outer automorphism is in the stabilizer of a primitive cohomology class in $H^1(R_n,\mathbb{Z}_k)$. In Section~\ref{sec:normal_series}, we derive a normal series for $\Out_{F_{n_k}}(F_n)$, and in Sections~\ref{sec:maxima} and~\ref{sec:self_norm}, we prove that $\Out_{F_{n_k}}(F_n)$ is maximal and self-normalizing in $\Out(F_n)$, respectively. In Section~\ref{sec:gen_set}, we determine a generating set for $\Out_{F_{n_k}}(F_n)$. Finally, in Section~\ref{sec:Inf_cover}, we consider the infinite-sheeted cover $q_n:\mathcal{C}\longrightarrow R_n$ and obtain a generating set for $\Out_{F_\infty}(F_n)$.

\section{Preliminaries}\label{sec:preliminaries}
In this section, we introduce several concepts related to $\Out(F_n)$ that are central to the theory developed in this paper.
\subsection{Topological representatives for automorphisms and outer automorphisms}
A continuous, basepoint-preserving map $f : (X,p) \to (X,p)$ is a \emph{self-homotopy equivalence} if it has a basepoint-preserving homotopy inverse $h : (X,p) \to (X,p)$ such that $f \circ h \sim \operatorname{id}_X \ \text{rel.\ } p\ \text{and}\ h \circ f \sim \operatorname{id}_X \ \text{rel.\ } p$.
The \emph{mapping class group} of $(X,p)$ is the group
\[
M(X,p) := \{\, [f]_p \mid f : X \to X \text{ is a homotopy equivalence
rel. } p \,\}.
\]
We now specialize to $X = R_n$, the wedge of $n$ circles, where $\pi_1(R_n)\cong F_n$, and the mapping class group of basepoint-preserving self-homotopy equivalences is naturally isomorphic to the automorphism group of the free group $F_n$, that is, $\Aut(F_n)\cong M(R_n,p)$. Similarly, when the basepoint is not fixed, the mapping class group is naturally isomorphic to the outer automorphism group of $F_n$, giving $\Out(F_n)\cong M(R_n)$ (see~\cite{CullerVogtmann}). Thus, passing from $M(R_n,p)$ to $M(R_n)$ corresponds to passing from $\Aut(F_n)$ to $\Out(F_n)$ by quotienting out the inner automorphisms.

\subsection{\texorpdfstring{A generating set for $\Aut(F_n)$}{A generating set for Aut(Fn)}}
Recall the notation from Section~\ref{sec:intro}. Algebraically, the automorphism group $\Aut(F_n)$ is generated by a finite collection of elementary Nielsen automorphisms. For distinct indices $1\leq i,j \leq n$, define the following automorphisms of $F_n$:
\[
\iota_i:
\begin{cases}
x_i\mapsto x_i^{-1},\\
x_k\mapsto x_k, & k\neq i,
\end{cases}
\qquad
\sigma_{ij}:
\begin{cases}
x_i\mapsto x_j,\\
x_j\mapsto x_i,\\
x_k\mapsto x_k, & k\notin\{i,j\},
\end{cases}
\]
and
\[
\mu_{ij}:
\begin{cases}
x_i\mapsto x_i x_j,\\
x_k\mapsto x_k, & k\neq i,
\end{cases}
\qquad
\nu_{ij}:
\begin{cases}
x_i\mapsto x_j x_i,\\
x_k\mapsto x_k, & k\neq i.
\end{cases}
\]
Here, $\iota_i$ is an inversion of the generator $x_i$, $\sigma_{ij}$ interchanges the generators $x_i$ and $x_j$, while $\mu_{ij}$ and $\nu_{ij}$ are, respectively, right and left Nielsen transformations. These elementary automorphisms generate $\Aut(F_n)$. We use the same notation for the images of these automorphisms under the natural projection $\operatorname{proj}_n:\Aut(F_n)\longrightarrow\Out(F_n)$. Thus, the corresponding elements of $\Out(F_n)$ will also be denoted by $\iota_i,\sigma_{ij},\mu_{ij}$, and $\nu_{ij}$, whenever there is no ambiguity. These generators provide a convenient finite generating set for both $\Aut(F_n)$ and, after passing to the quotient by inner automorphisms, $\Out(F_n)$. Gersten~\cite{Gersten1984} gave a finite presentation of $\Out(F_n)$ in terms of these generators.

\begin{theorem}[Gersten~\cite{Gersten1984}]
Suppose $n \geqslant 3$. The group $\Out(F_n)$ is generated by
\[
\{\iota_1, \mu_{ij}, \nu_{ij} \mid i,j = 1, \dots, n, \ i \neq j\},
\]
with relations
\begin{enumerate}[(i)]
\item $[\mu_{ij}, \mu_{kl}] = [\nu_{ij}, \nu_{kl}] = 1$ for
$k \notin \{i,j\}$, $l \neq i$;
\item $[\nu_{ij}, \mu_{kl}] = 1$ for $k \neq j$, $l \neq i$;
\item $[\mu_{ij}^{-1}, \mu_{jk}^{-1}] = [\mu_{ij}, \nu_{jk}]
      = [\mu_{ij}^{-1}, \mu_{jk}]^{-1}
      = [\mu_{ij}, \nu_{jk}^{-1}]^{-1} = \mu_{ik}^{-1}$
      for $k \notin \{i,j\}$;
\item $[\nu_{ij}^{-1}, \nu_{jk}^{-1}] = [\nu_{ij}, \mu_{jk}]
      = [\nu_{ij}^{-1}, \nu_{jk}]^{-1}
      = [\nu_{ij}, \mu_{jk}^{-1}]^{-1} = \nu_{ik}^{-1}$
      for $k \notin \{i,j\}$;
\item $\mu_{ij}\mu_{ji}^{-1}\nu_{ij} = \nu_{ij}\nu_{ji}^{-1}\mu_{ij}$,
      \quad $(\mu_{ij}\mu_{ji}^{-1}\nu_{ij})^4 = 1$;
\item $[\iota_1, \mu_{ij}] = [\iota_1, \nu_{ij}] = 1$ for $i,j \neq 1$;
\item $\mu_{12}^{\iota_1} = \nu_{12}^{-1}$, \quad
      $\mu_{21}^{\iota_1} = \mu_{21}^{-1}$;
\item $\iota_1^2 = 1$;
\item $\prod_{i \neq j} \mu_{ij}\nu_{ij}^{-1} = 1$ for each fixed $j$.
\end{enumerate}
\end{theorem}
With this presentation in hand, we return to the topological picture, where $\Out(F_n)$ is identified with $M(X)$, the group of free self-homotopy equivalences of a graph $X$ with $\pi_1(X)\cong F_n$,
taken up to free homotopy through self-homotopy equivalences. Two self-homotopy equivalences $g_1,g_2:X\to X$ represent the same element of $M(X)$ if and only if there is a free self-homotopy $H : g_1\sim g_2$ such that $x\mapsto H(x,t)$ is a free self-homotopy equivalence for every $t$. Let $p : X\to R_n$ be a normal covering map. This produces a short exact sequence
\begin{equation}
1\longrightarrow D\longrightarrow M_p^*(X)\longrightarrow M_p(R_n)\longrightarrow 1,
\end{equation}
where $M_p(R_n)\leq M(R_n)$ consists of those free homotopy equivalences $h : R_n \to R_n$ that lift to some $h' : X\to X$ with $p\circ h' = h\circ p$, up to homotopy; $M_p^*(X)$ is the group of all such lifts, up to homotopy; and $D$ is the deck transformation group of $p$ (see \cite{Sasao}).

Dowdall and Taylor~\cite{spence} gave an algebraic description of the short exact sequence associated with finite-sheeted unbranched covers, thereby establishing a Birman--Hilden analogue for $\Out(F_n)$ in the case of normal covers. Let $E\leq F_n$ be a finite-index subgroup, and let $\Gamma\leq\Out(F_n)$ be finitely generated. Let $\Gamma_E$ be the subgroup of $\Gamma$ consisting of all outer automorphisms that preserve the conjugacy class of $E$. Let $\Gamma^E$ be the induced subgroup of $\Out(E)$. That is, $f \in \Gamma^E$ if there exists $\phi \in \Aut(F_n)$ whose outer class lies in $\Gamma$ such that the restriction $\phi|_E$ represents the outer class $f$.
Dowdall and Taylor showed that these fit into a short exact sequence.
\begin{lemma}[Dowdall--Taylor~\cite{spence}]
\label{lem:short_exact_sequence}
Let $n\geq2$, let $E\leq F_n$ be a finite-index subgroup, and let $\Gamma\leq\Out(F_n)$ be finitely generated. Then the following sequence is exact:
\begin{equation}
\label{eq:BH-OutFn}
1 \longrightarrow N(E)/E\longrightarrow \Gamma^E\longrightarrow \Gamma_E\longrightarrow 1.
\end{equation}
Here $N(E)$ denotes the normalizer of $E$ in $F_n$, and the quotient
$N(E)/E$ is a finite group.
\end{lemma}
\subsection{Birman--Hilden analogue for regular cyclic covers}
Let $R_n$ be the rose with $n$ petals, so that $\pi_1(R_n)\cong F_n$. Let $p_k\colon \tilde R_n\to R_n$ be the regular covering induced by a free $\mathbb Z_k$-action. Then $\tilde R_n$ is a graph with $\pi_1(\tilde R_n)\cong F_{n_k}$, where $n_k=k(n-1)+1$. Here we identify $F_{n_k}$ with $E\leq F_n$ via $(p_k)_*$
and take $\Gamma=\Out(F_n)$. We now describe the three terms arising in the Birman--Hilden sequence.

\medskip
\noindent\textit{(i) The kernel $N(E)/E$.}
Since the cover is regular, $E\trianglelefteq F_n$, and hence $N(E)=F_n$. Thus,
\[N(E)/E\cong F_n/E\cong\mathbb Z_k,\]
which is precisely the \emph{deck transformation group} of the covering
$p_k$.

\medskip
\noindent\textit{(ii) The symmetric subgroup $\Gamma^E$.}
The symmetric subgroup consists of the outer automorphisms of $E$
induced by automorphisms of the ambient group $F_n$ that preserve $E$. Thus,
\[\Gamma^E=\Out^{F_n}(E)=p_E\bigl(\{\phi|_E:\phi\in\Aut(F_n),\ \phi(E)=E\}\bigr)\leq \Out(E),\]
where $p_E\colon\Aut(E)\to\Out(E)$ is the natural quotient map.

\medskip
\noindent\textit{(iii) The liftable subgroup $\Gamma_E$.}
The liftable subgroup consists of those outer automorphisms of $F_n$
that preserve the conjugacy class of $E$:
\[\Gamma_E=\Out_E(F_n)=\{[\phi]\in\Out(F_n):\phi(E)\text{ is conjugate to }E\}.\]
\medskip
With these identifications, the Birman--Hilden short exact sequence takes the following form.
\begin{cor}
\label{cor:birman_hilden_full_out}
Let $n\geq2$ and $k\geq2$, and let $p_k\colon\tilde R_n\to R_n$ be the regular $\mathbb Z_k$-cover
corresponding to $E\trianglelefteq F_n$, with $E\cong F_{n_k}$. Then there is a
short exact sequence
\[1\longrightarrow \mathbb Z_k\longrightarrow \Out^{F_n}(F_{n_k})
\longrightarrow \Out_{F_{n_k}}(F_n)\longrightarrow 1.\]
\end{cor}
\section{\texorpdfstring{Homological criteria for $\Out_{F_{n_k}}(F_n)$}{Homological criteria for the liftable subgroup}}\label{sec:homological_criteria}
Throughout, we use the notation introduced earlier. Let $e_1$ denote the first standard basis element of $H_1(R_n,\mathbb Z_k)\cong\mathbb Z_k^n$, and let $e_1^*$ denote its corresponding dual basis element in $H^1(R_n,\mathbb Z_k)\cong(\mathbb Z_k^n)^*$. Using basic covering space theory~\cite{Hatcher2002}, we obtain a homological characterization of $\Out_{F_{n_k}}(F_n)$ in the following theorem.
\begin{theorem}\label{thm:cyclic-cover-sequence}
A connected covering map $p_k:\widetilde R_n\to R_n$ is a regular $k$-sheeted cyclic cover if and only if there is a short exact sequence
\begin{equation}
1 \longrightarrow \pi_1\left(\widetilde{R}_n\right) \xrightarrow{(p_k)_*} \pi_1(R_n)\xrightarrow{\varphi}\mathbb{Z}_k\longrightarrow 1,
\end{equation}
where $\varphi$ is an epimorphism.
\end{theorem}
\subsection{\texorpdfstring{Criteria for $\Out_{F_{n_k}}(F_n)$}{Criteria for the liftable subgroup}}
Let $\phi \in \Out(F_n)$. Using the fact that our covering
action is cyclic, we obtain a homological criterion for the liftability
of $\phi$ with respect to the cover $p_k$.
\begin{theorem}\label{thm:homo_criteria}
Let $\phi\in\Out(F_n)$.
Then the following statements are equivalent.
\begin{enumerate}[(i)]
\item $\phi\in\Out_{F_{n_k}}(F_n)$.
\item If $\Psi(\phi)=(c_{ij})_{n\times n}$, then
$k\mid c_{1j}$ for $j\neq1$ and $\gcd(c_{11},k)=1$.
\item If $\Psi_k(\phi)=(e_{ij})_{n\times n}$, then
$e_{11}\in\mathbb Z_k^\times$ and $e_{1j}=0$ for $j\neq1$.
\item $\phi\in\Stab_{\Out(F_n)}
(\{\ell e_1^*:\ell\in\mathbb Z_k^\times\})$.
\end{enumerate}
\end{theorem}
\begin{proof}
$(i)\iff(ii)$  Let $\phi\in\Out(F_n)$. Let  $\varphi:\pi_1(R_n)\to \mathbb{Z}_k$
be an epimorphism such that $\ker(\varphi)=(p_k)_*(\pi_1(\widetilde R_n)).$
By Theorem~\ref{thm:cyclic-cover-sequence} and basic covering space theory, it follows that $\phi\in\Out_{F_{n_k}}(F_n) $  if and only if $\phi$  preserves the conjugacy class of  $(p_k)_*(\pi_1(\widetilde R_n))=\ker(\varphi)$.
Since $\mathbb Z_k$ is abelian, this is equivalent
to $\phi_\#(\ker(\bar\varphi))\subseteq
\ker(\bar\varphi)$, where $\bar{\varphi}$ is the induced homomorphism on the abelianized group $\pi_1(R_n)^{\mathrm{ab}} \cong H_1(R_n,\mathbb{Z})$. From the description of the induced map $(p_k)_{\#}$, we have
$\ker(\bar\varphi)=\langle ke_1,e_2,\ldots,e_n\rangle$, which proves our assertion.

The equivalence $(ii)\Longleftrightarrow(iii)$
follows immediately by reducing $\Psi(\phi)$
modulo $k$.
We now prove $(i)\Longleftrightarrow(iv)$.
Consider the following series of correspondences
between finite sets:
\begin{center}
$\displaystyle
\begin{gathered}
\mathcal C_k=
\left\{
\text{Regular }k\text{-sheeted covers }
p_k:\widetilde R_n\to R_n
\text{ whose deck transformation group is }
\mathbb Z_k
\right\}
\\[6pt]
\updownarrow
\\[6pt]
\left\{
\text{Monomorphisms}
(p_k)_*:\pi_1(\widetilde R_n)\hookrightarrow\pi_1(R_n)
\text{ such that }
\pi_1(R_n)/(p_k)_*(\pi_1(\widetilde R_n))
\cong\mathbb Z_k
\right\}
\\[6pt]
\updownarrow
\\[6pt]
\left\{
\left\{
\text{Epimorphisms }
\varphi:\pi_1(R_n)\to\mathbb Z_k:
(p_k)_*(\pi_1(\widetilde R_n))
=\ker(\varphi)
\right\}:p_k\in\mathcal C_k
\right\}
\\[6pt]
\updownarrow
\\[6pt]
\left\{
\left\{
\text{Epimorphisms }
\bar\varphi:H_1(R_n,\mathbb Z)\to\mathbb Z_k:
(p_k)_\#(H_1(\widetilde R_n,\mathbb Z))
=\ker(\bar\varphi)
\right\}:p_k\in\mathcal C_k
\right\}
\\[6pt]
\updownarrow
\\[6pt]
\left\{
\{\ell\alpha:\ell\in\mathbb Z_k^\times\}:
\alpha\in H^1(R_n,\mathbb Z_k)
\text{ is a primitive cohomology class}
\right\}.
\end{gathered}
$
\end{center}
$\Out(F_n)$ acts on $H^1(R_n,\mathbb Z_k)$ from the right via $\alpha\cdot\phi=\alpha\Psi_k(\phi)$. Since the above correspondences are $\Out(F_n)$-equivariant with respect to induced actions, the stabilizers of corresponding elements coincide. Thus, for the covering map $p_k$, which induces $(p_k)_\#$ as described in Equation~\eqref{eq:finite-cover-homology}, we obtain
\begin{equation}\label{eq:out_stab}
\Out_{F_{n_k}}(F_n)
=\Stab_{\Out(F_n)}(p_k)
=\Stab_{\Out(F_n)}
(\{\ell e_1^*:\ell\in\mathbb Z_k^\times\}),
\end{equation}
which proves $(i)\Longleftrightarrow(iv)$ and completes the proof.
\end{proof}
\begin{definition}\label{def:level-k-subgroup}
Let $\mathcal J_{n,k}=\rho_k(\operatorname{GL}(n,\mathbb Z))$,
where $\rho_k$ denotes reduction modulo $k$.
Consider the surjective homomorphism
$\Psi_k=\rho_k\circ\Psi:
\Out(F_n)\twoheadrightarrow
\mathcal J_{n,k}$. For $k \geq 2$, $\mathcal{T}_{n,k}:=\ker(\Psi_k)$.
\end{definition}
\noindent The following immediate consequence of Theorem~\ref{thm:homo_criteria} gives $\mathcal{T}_{n,k} \leq \Out_{F_{n_k}}(F_n)$.
\begin{cor}\label{cor:k-level-interse}
For $n \geq 2$ and $k, \ell \geq 2$, we have
\begin{enumerate}[(i)]
    \item $\mathcal{T}_{n,k}\triangleleft\Out_{F_{n_k}}(F_n)$.
    \item $\Out_{F_{n_k}}(F_n)\cap\Out_{F_{n_\ell}}(F_n)
    =\Out_{F_{n_d}}(F_n)$, where $d=\operatorname{lcm}(k,\ell)$.
\end{enumerate}
\end{cor}
Let $E$ denote Euler's totient function. Combining Theorem~\ref{thm:homo_criteria} with Corollary~\ref{cor:k-level-interse}, we obtain the following result.
\begin{cor}\label{cor:index_subgroup}
Let $n\geq 2$ and $k\geq 2$. Let $\mathcal{I}_{n,k}^{\mathbb Z}=\Psi(\Out_{F_{n_k}}(F_n))$ and
$\mathcal I_{n,k}=\Psi_k(\Out_{F_{n_k}}(F_n))$. Let $\mathcal J_{n,k}=\rho_k(\mathrm{GL}(n,\mathbb Z))$, where $\rho_k$ denotes reduction modulo $k$. Then the following statements hold.
\begin{enumerate}[(i)]
\item
$[\Out(F_n):\Out_{F_{n_k}}(F_n)]
=[\mathrm{GL}(n,\mathbb{Z}):\mathcal{I}_{n,k}^{\mathbb Z}]
=[\mathcal{J}_{n,k}:\mathcal{I}_{n,k}]$.
\item We have
\[
[\Out(F_n):\Out_{F_{n_k}}(F_n)]
=\frac{|\{\text{primitive elements of }(\mathbb Z_k^n)^*\}|}{E(k)}.
\]
In particular, if $k$ is prime, then
\[
[\Out(F_n):\Out_{F_{n_k}}(F_n)]=\frac{k^n-1}{k-1}.
\]
\end{enumerate}
\end{cor}
\subsection{\texorpdfstring{Criterion for $\Out_{F_{n_k}}(F_n,e_1^*)$}{Criterion for the stabilizer subgroup}}\label{sec:criter_stab}
Recall that $\Out_{F_{n_k}}(F_n,e_1^*)=\Stab_{\Out(F_n)}(e_1^*)$. We have the following criterion for membership in this subgroup.
\begin{theorem}\label{thm:crite_stab_e_1}
Let $\phi\in\Out(F_n)$. The following
statements are equivalent.
\begin{enumerate}[(i)]
\item $\phi\in \Out_{F_{n_k}}(F_n,e_1^*)$.
\item If $\Psi(\phi)=(c_{ij})_{n\times n}$,
then $k\mid c_{1j}$ for $1\leq j\leq n$
with $j\neq 1$, and $c_{11}\equiv 1\pmod{k}$.
\item If $\Psi_k(\phi)=(e_{ij})_{n\times n}$,
then $e_{11}=1$ and $e_{1j}=0$ for $1\leq j\leq n$ with $j\neq 1$.
\end{enumerate}
\end{theorem}
\begin{proof}
Let $A=(c_{ij})_{n\times n}=\Psi(\phi)
\in\mathrm{GL}(n,\mathbb Z)$ be the matrix associated to $\phi$, and let $A'=(e_{ij})_{n\times n}$ denote the reduction of $A$ modulo $k$, so that $A'=\Psi_k(\phi)\in\mathcal{J}_{n,k}$.

Suppose that $\phi\in \Out_{F_{n_k}}(F_n,e_1^*)$. Then by definition, $e_1^*\Psi_k(\phi)=e_1^*$. Since $e_1^*=(1,0,\ldots,0)$, this implies that the first row of $A'$ is $(1,0,\ldots,0)$. Equivalently,
$(e_{11},e_{12},\ldots,e_{1n})=(1,0,\ldots,0)$.
Hence, $e_{11}=1$ and $e_{1j}=0$ for all $j\neq 1$,
so $A'$ has the form
\[
A'=
\begin{bmatrix}
1 & 0 & 0 & \cdots & 0\\
e_{21} & e_{22} & e_{23} & \cdots & e_{2n}\\
\vdots & \vdots & \vdots & \ddots & \vdots\\
e_{n1} & e_{n2} & e_{n3} & \cdots & e_{nn}
\end{bmatrix}.
\]
Conversely, a matrix with this first row fixes $e_1^*$, so this proves $(i)\Longleftrightarrow(iii)$.
The equivalence $\text{(ii)}\Longleftrightarrow\text{(iii)}$
is immediate from reduction modulo $k$.
\end{proof}
\noindent For the proof of the next corollary, we will require the following lemma from basic group theory.
\begin{lemma}\label{lem:normality-under-surjection}
Let $\psi:G\to H$ be a surjective group
homomorphism and let $N\leq G$ be a subgroup
such that $N=\psi^{-1}(\psi(N))$. Then
\[
N\triangleleft G
\quad\Longleftrightarrow\quad
\psi(N)\triangleleft H.
\]
\end{lemma}
\begin{cor}\label{cor:index_stab_e1}
Let $n \ge 2$ and $k \ge 2$. Then
\begin{enumerate}
\item $\mathcal{T}_{n,k} \triangleleft \Out_{F_{n_k}}(F_n,e_1^*) \triangleleft \Out_{F_{n_k}}(F_n) $.
\item
$
\Out_{F_{n_k}}(F_n) /\Out_{F_{n_k}}(F_n,e_1^*) \cong \mathbb{Z}_k^{\times}.
$
In particular,
\[
[\Out_{F_{n_k}}(F_n):\Out_{F_{n_k}}(F_n,e_1^*)] = E(k).
\]
\end{enumerate}
\end{cor}
\begin{proof}
By Theorem~\ref{thm:crite_stab_e_1}, $\mathcal{T}_{n,k}\subset \Out_{F_{n_k}}(F_n,e_1^*)$ for every $k \ge 2$. From the definition of $\Out_{F_{n_k}}(F_n)$, it is also clear that $\Out_{F_{n_k}}(F_n,e_1^*)\le \Out_{F_{n_k}}(F_n)$. It remains to check that $\Out_{F_{n_k}}(F_n,e_1^*)$ is normal in $\Out_{F_{n_k}}(F_n)$. By Lemma~\ref{lem:normality-under-surjection}, this is equivalent to showing that
$\Psi_k(\Out_{F_{n_k}}(F_n,e_1^*)) \triangleleft \Psi_k(\Out_{F_{n_k}}(F_n))$.

Let $\phi \in \Out_{F_{n_k}}(F_n)$. By Theorem~\ref{thm:homo_criteria}, its image under $\Psi_k$ has the form $B=(e_{ij})_{n\times n}$ with $e_{1i}=0$ for $i\neq1$ and $e_{11}\in\mathbb Z_k^\times$. Consider the homomorphism
$\zeta : \Psi_k(\Out_{F_{n_k}}(F_n)) \longrightarrow \mathbb{Z}_k^{\times},\qquad\zeta(B)=e_{11}$. Then $\ker(\zeta)=\Psi_k(\Out_{F_{n_k}}(F_n,e_1^*))$.

Hence, $\Psi_k(\Out_{F_{n_k}}(F_n,e_1^*))$ is a normal subgroup of $\Psi_k(\Out_{F_{n_k}}(F_n))$, with quotient isomorphic to $\mathbb{Z}_k^{\times}$. The second and first isomorphism theorems then imply that $\Out_{F_{n_k}}(F_n,e_1^*)$ is a normal subgroup of $\Out_{F_{n_k}}(F_n)$ and that $\Out_{F_{n_k}}(F_n)/\Out_{F_{n_k}}(F_n,e_1^*) \cong \mathbb{Z}_k^{\times}$.
Since $|\mathbb{Z}_k^{\times}|=E(k)$, the index formula follows.
\end{proof}
\section{Applications of the homology representation}\label{sec:applic}
We apply Theorems~\ref{thm:homo_criteria} and~\ref{thm:crite_stab_e_1} to obtain a normal series for the liftable subgroup. We then prove self-normalization for $n\geq3$ and characterize maximality for $n\geq2$.
\subsection{\texorpdfstring{Normal series of subgroups of $\Out(F_n)$}{Normal series of subgroups of Out(Fn)}}\label{sec:normal_series}
Putting together the results in Corollaries~\ref{cor:index_subgroup} and~\ref{cor:index_stab_e1}, we obtain the following normal series of subgroups of $\Out(F_n)$.
\begin{theorem}\label{thm:normal_series_app}
Let $n\ge2$ and $k\ge2$. Then  $\Out_{F_{n_k}}(F_n)$ admits a normal series
\[ 1 \triangleleft \mathcal{T}_{n,k} \triangleleft \Out_{F_{n_k}}(F_n,e_1^*) \triangleleft \Out_{F_{n_k}}(F_n).\] Furthermore,
\[\Out_{F_{n_k}}(F_n) /\Out_{F_{n_k}}(F_n,e_1^*) \cong \mathbb{Z}_k^{\times}\]
and
\[
[\Out(F_n):\Out_{F_{n_k}}(F_n,e_1^*)]
=
\left|\left\{\text{primitive elements of }(\mathbb{Z}_k^n)^*\right\}\right|
.
\]
\end{theorem}
It follows from Theorem~\ref{thm:crite_stab_e_1} and Corollary~\ref{cor:index_stab_e1}, by the first and third isomorphism theorems, that the surjective homomorphism $\Psi_k$ induces an isomorphism
\[
\overline\Psi_k:
\frac{\Out_{F_{n_k}}(F_n)}{\Out_{F_{n_k}}(F_n,e_1^*)}
\longrightarrow
\frac{\Psi_k(\Out_{F_{n_k}}(F_n))}{\Psi_k(\Out_{F_{n_k}}(F_n,e_1^*))}.
\]
Both quotient groups are isomorphic to the multiplicative group $\mathbb{Z}_k^\times$. Let $\ell\in\mathbb{Z}_k^\times$, and denote its multiplicative inverse by $\bar\ell$. Fix integer representatives of these residue classes whenever they occur as exponents. Consider the matrix $D_\ell = ((d_\ell)_{ij})_{n\times n} \in \mathrm{GL}(n,\mathbb{Z}_k)$, where
\[
(d_\ell)_{ij}=
\begin{cases}
\ell, & i=j=1,\\
\bar{\ell}, & i=j=2,\\
1, & i=j\geq 3,\\
0, & i\neq j.
\end{cases}
\]

Note that $D_\ell$ represents a coset in $\Psi_k(\Out_{F_{n_k}}(F_n))/\Psi_k(\Out_{F_{n_k}}(F_n,e_1^*))$. Moreover, each coset represented by a $D_\ell$ pulls back under $\overline{\Psi}_k$ to a coset represented by $w_\ell=\mu_{12}^{1-\bar{\ell}}\mu_{21}^{-1}\mu_{12}^{1-\ell}\mu_{21}^{\bar{\ell}}$.
Hence, $w_\ell$ represents the corresponding coset in
$\Out_{F_{n_k}}(F_n)/
\Out_{F_{n_k}}(F_n,e_1^*)$.
\begin{cor}\label{cor:coset-representatives}
For $n\geq 2$ and $k\geq 2$, the distinct cosets of
$\Out_{F_{n_k}}(F_n)/
\Out_{F_{n_k}}(F_n,e_1^*)$ are represented by the elements
\[
\mathcal{S}_k'' := \left\{
w_\ell = \mu_{12}^{1-\bar{\ell}}\mu_{21}^{-1}
\mu_{12}^{1-\ell}\mu_{21}^{\bar{\ell}}
: \ell\in\mathbb{Z}_k^\times
\right\}.
\]
\end{cor}
\subsection{An involution and the normal series}\label{sec:normal_involu}
Let $n\geq2$, and recall the basis-inverting involution $\iota=\iota_1\cdots\iota_n\in\Out(F_n)$. Let $I_n\in\mathcal J_{n,k}$ denote the identity matrix. Since $\Psi_k(\iota) = -I_{n}$, it follows from Theorem~\ref{thm:crite_stab_e_1} that $\iota \in \Out_{F_{n_k}}(F_n,e_1^*)$ if and only if $k = 2$. Moreover, for $k = 2$, it is apparent from Theorem~\ref{thm:normal_series_app} that
$\Out_{F_{n_k}}(F_n)
= \Out_{F_{n_k}}(F_n,e_1^*)$. Thus, as $\Psi_k(\iota)$
commutes with every element of $\mathcal{J}_{n,k}$,
we obtain the following extension of the normal series in Theorem~\ref{thm:normal_series_app}.
\begin{cor}\label{cor:norm_invo}
For $n \geq 2$ and $k \geq 3$, the group $\Out_{F_{n_k}}(F_n)$ admits the normal series
\[
1
\triangleleft \mathcal{T}_{n,k}
\triangleleft \Out_{F_{n_k}}(F_n,e_1^*)
\triangleleft \langle\!\langle \Out_{F_{n_k}}(F_n,e_1^*),\iota
\rangle\!\rangle
\triangleleft \Out_{F_{n_k}}(F_n),
\]
where
$\langle\!\langle \Out_{F_{n_k}}(F_n,e_1^*),\iota
\rangle\!\rangle$ denotes the normal closure of
$\langle \Out_{F_{n_k}}(F_n,e_1^*),\iota\rangle$ in
$\Out_{F_{n_k}}(F_n)$,
\[
\langle\!\langle \Out_{F_{n_k}}(F_n,e_1^*),\iota
\rangle\!\rangle /
\Out_{F_{n_k}}(F_n,e_1^*)
\cong \mathbb{Z}_2,
\]
and
\[
\left[
\Out_{F_{n_k}}(F_n):
\langle\!\langle \Out_{F_{n_k}}(F_n,e_1^*),\iota
\rangle\!\rangle
\right]
= \frac{E(k)}{2}.
\]
\end{cor}
\noindent From Corollary~\ref{cor:norm_invo}, it follows that
$\Out_{F_{n_k}}(F_n)
= \langle\!\langle
\Out_{F_{n_k}}(F_n,e_1^*),\iota
\rangle\!\rangle$
if and only if $E(k)=2$. Hence, we have the following corollary.
\begin{cor}\label{cor:i_in_stab_e1}
Let $n\geq2$ and $k\geq3$. Then
$\Out_{F_{n_k}}(F_n)
= \langle\!\langle \Out_{F_{n_k}}(F_n,e_1^*),\iota
\rangle\!\rangle$
if and only if $k \in \{3,4,6\}$.
\end{cor}
\subsection{\texorpdfstring{Maximality of $\Out_{F_{n_k}}(F_n)$}{Maximality of the liftable subgroup}}\label{sec:maxima}
To establish the maximality of $\Out_{F_{n_k}}(F_n)$, we recall some basic notions from permutation group theory. We begin with transitive group actions, blocks, and primitive actions, and then introduce the projective space $\mathbb{P}^{n-1}(\mathbb{Z}_k)$ over the finite field $\mathbb{Z}_k$. We take the right action of $\Out(F_n)$ on this space and recall a standard characterization of primitive actions in terms of maximal point stabilizers. These results will be used to derive a group-theoretic corollary needed in the proof of our main theorem.
\begin{definition}\label{def:group_th_max}
Let $G$ be a group acting on a set $\Delta$. We recall the following standard terminology.
\begin{enumerate}[(i)]
\item Suppose that the action is transitive. A nonempty subset $B\subseteq\Delta$ is called a \emph{block} if, for every $g\in G$, either $(B)\cdot g=B$ or $(B)\cdot g\cap B=\varnothing$.
\item The blocks $\{x\}$, for $x\in\Delta$, and $\Delta$ itself are called \emph{trivial blocks}. All other blocks are called \emph{nontrivial blocks}.
\item A transitive action of $G$ on $\Delta$ is called \emph{primitive} if it admits no nontrivial blocks.
\item For $x\in\Delta$, the subgroup $\Stab_{G}(x)=\{g\in G\mid x\cdot g=x\}$ is called the \emph{point stabilizer} of $x$ in $G$.
\end{enumerate}
\end{definition}

\noindent We next recall the notion of projective space over a finite field.
\begin{definition}\label{def:projective_space}
Let $k$ be a prime and let $\mathbb{Z}_k$ denote the finite field with $k$ elements. The \emph{projective space of dimension $n-1$ over $\mathbb{Z}_k$}, denoted by $\mathbb{P}^{n-1}(\mathbb{Z}_k)$, is defined as
\[
\mathbb{P}^{n-1}(\mathbb{Z}_k)
=
((\mathbb Z_k^n)^*\setminus\{0\})/\sim,
\]
where $v\sim w$ if and only if $v=\lambda w$ for some $\lambda\in\mathbb{Z}_k^\times$. Each equivalence class $[\beta]$ corresponds to the one-dimensional subspace spanned by $\beta$.
\end{definition}
\noindent The following lemma records the transitivity of the projective action.
\begin{lemma}\label{lem:Gl_n_act_tran_proj}
The group $\Out(F_n)$ acts on $\mathbb{P}^{n-1}(\mathbb{Z}_k)$ by $ \beta \cdot \phi = \beta \Psi_k(\phi)$, and this action is transitive.
\end{lemma}

\begin{cor}[Corollary 1.5(A) of \cite{dixon_mortimer_1996}]
\label{cor:primitive_maximal_stabilizer}
Let $G$ act transitively on a set $\Delta$ with at least two elements. Then the action of $G$ on $\Delta$ is primitive if and only if the point stabilizer $\Stab_{G}(x)$ is a maximal subgroup of $G$, for some (equivalently, every) $x\in\Delta$.
\end{cor}
\begin{lemma}\label{lem:max_point_stab}
Let $k$ be prime and $n\geq2$. For every line $L \in \mathbb{P}^{n-1}(\mathbb{Z}_k)$, the point stabilizer $\Stab_{\Out(F_n)}(L)$ is a maximal subgroup of $\Out(F_n)$.
\end{lemma}
\begin{proof}
By Corollary~\ref{cor:primitive_maximal_stabilizer}, it suffices to prove that the action of $\Out(F_n)$ on $\mathbb P^{n-1}(\mathbb Z_k)$ is primitive. First, we show that the stabilizer of a point $L_0$ acts transitively on its complement. Write $L_0=\langle e\rangle$, and let $L_1,L_2\neq L_0$ have nonzero representatives $v,v'$, respectively. Extend $(e,v)$ and $(e,v')$ to ordered bases of $(\mathbb Z_k^n)^*$. Let $A\in\mathrm{GL}(n,\mathbb Z_k)$ be the matrix taking the first basis to the second under the right action, and put $d=\det A$. Replace the second vector $v'$ in the target basis by $d^{-1}v'$, leaving the other target basis vectors unchanged. The resulting matrix $A_0$ has determinant $1$, fixes $e$, and maps $v$ to $d^{-1}v'$. This construction also applies when $n=2$, since scaling $v'$ does not change $L_2$.

Every matrix in $\mathrm{SL}(n,\mathbb Z_k)$ is a product of elementary matrices. Each elementary matrix $I+aE_{ij}$ lifts to the integral matrix $I+\widetilde aE_{ij}$, where $\widetilde a$ represents $a$, and hence to a power of a Nielsen transformation. Thus $A_0\in\mathcal J_{n,k}$, and there exists $\phi\in\Out(F_n)$ with $\Psi_k(\phi)=A_0$. This outer automorphism fixes $L_0$ and sends $L_1$ to $L_2$. Therefore $\Stab_{\Out(F_n)}(L_0)$ acts transitively on $\mathbb P^{n-1}(\mathbb Z_k)\setminus\{L_0\}$.

Let $B\subseteq\mathbb P^{n-1}(\mathbb Z_k)$ be a block containing $L_0$. If $B$ contains no lines other than $L_0$, then $B = \{L_0\}$, which is a trivial block. If $B \neq \{L_0\}$, there exists a line $L_1 \in B$ with $L_1 \neq L_0$. For any line $L_2 \in \mathbb{P}^{n-1}(\mathbb{Z}_k)\setminus\{L_0\}$, the transitivity of $\Stab_{\Out(F_n)}(L_0)$ on $\mathbb{P}^{n-1}(\mathbb{Z}_k) \setminus \{L_0\}$ gives an element $g\in\Stab_{\Out(F_n)}(L_0)$ such that $(L_1)\cdot g = L_2$. Since $g \in \Stab_{\Out(F_n)}(L_0) $ , $(L_0)\cdot g = L_0 \in B$, which implies $L_0 \in B \cap (B)\cdot g$ and thus $B \cap (B)\cdot g \neq \emptyset$.
By the definition of a block, this nonempty intersection implies $(B)\cdot g = B$. Since $L_1 \in B$, it follows that $L_2 = (L_1)\cdot g \in (B)\cdot g = B$. Since $L_2$ was arbitrarily chosen from $\mathbb{P}^{n-1}(\mathbb{Z}_k) \setminus \{L_0\}$, we have $\mathbb{P}^{n-1}(\mathbb{Z}_k) \setminus \{L_0\} \subseteq B$, which gives $B = \mathbb{P}^{n-1}(\mathbb{Z}_k)$, which is a trivial block. Since all blocks containing $L_0$ are trivial, $\mathbb{P}^{n-1}(\mathbb{Z}_k)$ contains no nontrivial blocks, and the action of $\Out(F_n)$ on $\mathbb{P}^{n-1}(\mathbb{Z}_k)$ is primitive. By Corollary~\ref{cor:primitive_maximal_stabilizer}, $\Stab_{\Out(F_n)}(L)$ is a maximal subgroup of $ \Out(F_n)$ for every $L\in\mathbb{P}^{n-1}(\mathbb{Z}_k)$.
\end{proof}

We can now prove the maximality theorem.

\begin{theorem}\label{thm:maximal}
For $n\ge 2$, the subgroup $\Out_{F_{n_k}}(F_n)$ is maximal in $\Out(F_n)$ if and only if $k$ is prime.
\end{theorem}

\begin{proof}
If $k$ is prime, then Theorem~\ref{thm:homo_criteria} (iv) gives
$\Out_{F_{n_k}}(F_n)=\Stab_{\Out(F_n)}(\{\ell e_1^* : \ell \in \mathbb{Z}_k^\times\})$.
Let $L_1:=\operatorname{span}_{\mathbb{Z}_k}(e_1^*)\in\mathbb{P}^{n-1}(\mathbb{Z}_k)$. Since
$\{\ell e_1^* : \ell \in \mathbb{Z}_k^\times\}=L_1\setminus\{0\}$, we obtain
$\Out_{F_{n_k}}(F_n)=\Stab_{\Out(F_n)}(L_1)$.
Thus, $\Out_{F_{n_k}}(F_n)$ is the stabilizer of the line $L_1$, and the desired assertion follows from Lemma~\ref{lem:max_point_stab}.

Conversely, suppose that $k$ is composite, and let $p$ be a prime divisor of $k$. The homological criterion gives $H_k\leq H_p$, where $H_j=\Out_{F_{n_j}}(F_n)$. The automorphism $\mu_{21}^{p}$ has homology matrix $I+pE_{12}$, so it belongs to $H_p$ but not to $H_k$, since $k\nmid p$. Moreover, $\mu_{21}\notin H_p$. Hence $H_k<H_p<\Out(F_n)$, proving that $H_k$ is not maximal.
\end{proof}
\subsection{\texorpdfstring{$\Out_{F_{n_k}}(F_n)$ is self-normalizing}{The liftable subgroup is self-normalizing}}\label{sec:self_norm}
A subgroup $H$ of a group $G$ is said to be \textit{self-normalizing} if $N_G(H)=H$. In other words, this means that if $H\trianglelefteq K\leq G$, then $H=K$. In this subsection, we show that $
\Out_{F_{n_k}}(F_n)$ satisfies this property.

\begin{theorem}\label{thm:self_norm}
For $n\geq3$ and $k\geq2$, $\Out_{F_{n_k}}(F_n)$ is self-normalizing in $\Out(F_n)$.
\end{theorem}

\begin{proof}
Suppose, to the contrary, that
$N_{\Out(F_n)}(\Out_{F_{n_k}}(F_n))\neq \Out_{F_{n_k}}(F_n)$,
and choose
$h\in N_{\Out(F_n)}(\Out_{F_{n_k}}(F_n))
\setminus \Out_{F_{n_k}}(F_n)$.
By the definition of the normalizer,
$h^{-1}\phi h\in\Out_{F_{n_k}}(F_n)$
for every $\phi\in\Out_{F_{n_k}}(F_n)$.
By Theorem~\ref{thm:homo_criteria}, this is equivalent to requiring that
$e_1^*\cdot(h^{-1}\phi h)=\ell e_1^*$
for some $\ell\in\mathbb{Z}_k^\times$. Equivalently,
$(e_1^*\cdot h^{-1})\cdot\phi=\ell(e_1^*\cdot h^{-1})$.
Setting $\beta:=e_1^*\cdot h^{-1}\in(\mathbb{Z}_k^n)^*$, we obtain
$\beta\cdot\phi=\ell\beta$
for some $\ell\in\mathbb{Z}_k^\times$. Thus, it suffices to show that for every primitive cohomology class $\beta\in H^1(R_n,\mathbb{Z}_k)$ satisfying
$\beta\notin\{\ell e_1^*:\ell\in\mathbb{Z}_k^\times\}$, there exists $\phi\in\Out_{F_{n_k}}(F_n)$ such that
$\beta\cdot\phi\neq\ell\beta$
for every $\ell\in\mathbb{Z}_k^\times$.
We get a contradiction by choosing an elementary matrix that arises as the image under $\Psi$ of a Nielsen transformation in $\Out(F_n)$.

We first show that $\beta_1$ is a unit in $\mathbb Z_k$,
where $\beta_1$ denotes the first entry of the cohomology class
$\beta$. Suppose $\beta_1$ is not a unit. Then there exists a prime $p\mid k$ such that $p\mid\beta_1$. Since $\beta$ is primitive, there exists
some $i\in\{2,\ldots,n\}$ such that $p\nmid\beta_i$. Now take $\phi=\mu_{1i}\in\Out_{F_{n_k}}(F_n)$. Then \[\beta\cdot\phi=(\beta_1+\beta_i,\beta_2,\ldots,\beta_n)=\ell\beta\pmod{k}\] for some $\ell\in\mathbb Z_k^\times$.
Equating the first coordinates on both sides, we obtain $\beta_1+\beta_i\equiv\ell\beta_1\pmod{k}$, which implies
$\beta_i\equiv(\ell-1)\beta_1\pmod{k}$. Since $p\mid\beta_1$, it follows that $p\mid\beta_i$, which contradicts our choice of $i$. Consequently, $\beta_1$ must be a unit in $\mathbb Z_k$.

Now, since $\beta\notin\{\ell e_1^*:\ell\in\mathbb Z_k^\times\}$,
there exists some $i\geq2$ such that $\beta_i\neq0$ in $\mathbb Z_k$.
Since $n\geq3$, we choose $j\geq2$ with $j\neq i$. Take $\phi=\mu_{ji}\in\Out_{F_{n_k}}(F_n)$. Then \[\beta\cdot\phi=(\beta_1,\ldots,\beta_j+\beta_i,\ldots,\beta_n)=\ell\beta\pmod{k}\] for some $\ell\in\mathbb Z_k^\times$.
Equating the first coordinates on both sides gives $\beta_1=\ell\beta_1\pmod{k}$. Since $\beta_1$ is a unit, we obtain $\ell=1$. Equating the $j$-th coordinates on both sides now gives $\beta_j+\beta_i=\beta_j\pmod{k}$, which implies $\beta_i=0$ in $\mathbb Z_k$, which contradicts our choice of $i$.

This contradiction completes the proof.
\end{proof}
\subsection{Generating the homological image of the stabilizer}\label{sec:gen_set}
We show the following theorem in this subsection.
\begin{theorem}\label{thm:gen_img_stab_e1}
Let $n\geq3$ and $k\geq2$. The group $\Psi(\Out_{F_{n_k}}(F_n,e_1^*))$ is generated by
\[\mathcal S^{(n)} = \kappa_n\big(\Gamma(k)\big) \,\cup\, \Psi(S),\]
where $\kappa_n : \mathrm{SL}(2,\mathbb{Z}) \to \mathrm{GL}(n,\mathbb{Z})$ denotes the canonical embedding defined by
\[A \mapsto \begin{bmatrix}A & O_{2,\,n-2} \\O_{n-2,\,2} & I_{n-2}\end{bmatrix}.\]
Here, $O_{m,n}$ denotes the $m\times n$ zero matrix, and $I_{n-2}$
denotes the identity matrix of size $(n-2)\times(n-2)$.
The generating set $S$ appearing in the statement of the theorem is
given explicitly by
\begin{equation}\label{eq:stabilizer-image-generators}
S=\{\mu_{ij} \mid 2 \le i,j \le n,\ i \neq j\}
\cup \{\sigma_{ij}\mid2\leq i<j\leq n\}
\cup \{\iota_2,\mu_{12}\}.
\end{equation}
\end{theorem}
\begin{rem}
We use the convention $\mu_{ij}\colon x_i\mapsto x_ix_j$ and
$\Psi(\mu_{ij})=I+E_{ji}$, where $E_{ji}$ is the matrix with
a $1$ in the $(j,i)$-entry and zeros elsewhere. Thus, right
multiplication by $\Psi(\mu_{ij})$ performs $C_i\to C_i+C_j$,
while left multiplication performs $R_j\to R_j+R_i$.
\end{rem}

The proof of Theorem~\ref{thm:gen_img_stab_e1} relies on the following lemma.

\begin{lemma}\label{lem:subrow-reduction}
Let $n \ge 3$, and let $A = (a_{ij}) \in \mathrm{GL}(n, \mathbb{Z})$. Let $\mathbf{a} = (a_{12}, a_{13}, \dots, a_{1n}) \in \mathbb{Z}^{n-1}$ be the subrow vector consisting of entries from column $2$ through $n$ in row $1$, and let $L = \gcd(a_{12}, \dots, a_{1n}) \in \mathbb{Z}_{\ge0}$.
Define the generator subset $S' \subset S \subset \Out(F_n)$ by $ S' = \{\mu_{ij} \mid 2 \le i, j \le n, i \neq j\} \cup \{\sigma_{ij}\mid2\leq i<j\leq n\} \cup \{\iota_{2}\}$. Then there exists a word $w\in\langle S'\rangle$ such that $g=\Psi(w)\in\mathrm{GL}(n,\mathbb Z)$ satisfies $(Ag)_{12} = L$, and $(Ag)_{1i} = 0$ for $3\le i \le n$.
\end{lemma}
\begin{proof}
We divide the argument into two cases.

\textit{Case 1: Degenerate case.} Suppose that $a_{12}=a_{13}=\cdots=a_{1n}=0$. Then $L=0$, so the empty word $w=1$ and $g=I_n$ satisfy the conclusion.

\textit{Case 2: Nondegenerate case.} Suppose that $a_{12},a_{13},\ldots,a_{1n}$ are not all zero. We proceed by induction on $m \in \{3, \dots, n\}$. For each $m$, let $P(m)$ be the statement that there exists a word $w_m \in \langle S' \rangle$ such that the matrix $A^{(m)} := A\Psi(w_m) = (a_{ij}^{(m)})$ satisfies
\[
a_{1i}^{(m)}=0\quad(3\leq i\leq m),
\qquad |a_{12}^{(m)}|=\gcd(a_{12},\ldots,a_{1m}).
\]
We construct $w_m$ using only column operations among columns $2,\ldots,m$, so column $1$ and all columns beyond $m$ remain unchanged.

If $a_{13} = 0$, then $P(3)$ holds trivially with the empty word $w_3 = \mathrm{id}$. If $a_{13} \neq 0$, applying the Euclidean algorithm to columns $2$ and $3$ via right multiplication by the images of the Nielsen transformation $\mu_{23}$, its inverse, and the transposition $\sigma_{23}$ changes the pair $(a_{12}, a_{13})$ into $(\pm\gcd(a_{12}, a_{13}), 0)$. Thus there exists $w_3\in\langle S'\rangle$ satisfying both conditions in $P(3)$.

Now assume $P(m)$ holds for some $3 \le m < n$, so that $a_{13}^{(m)}=\cdots=a_{1m}^{(m)}=0$. We eliminate the entry $a_{1(m+1)}^{(m)}$ using $a_{12}^{(m)}$ via the Euclidean algorithm. Let $x = a_{12}^{(m)}$ and $y = a_{1(m+1)}^{(m)}$. If $y = 0$, then $P(m+1)$ holds immediately with $w_{m+1} = w_m$. If $y \neq 0$, we apply the following steps iteratively:
\begin{enumerate}[(i)]
\item Apply the division algorithm to write $x = q \cdot y + r$, where $q, r \in \mathbb{Z}$ and $0 \le r < |y|$.
\item Right-multiply by $\Psi(\mu_{2(m+1)})^{-q}$ to perform $C_2\to C_2-qC_{m+1}$, replacing the $(1,2)$ entry with $r$.
\item Interchange $C_2$ and $C_{m+1}$ by right-multiplying by $\Psi(\sigma_{2(m+1)})$, changing $(r,y)$ into $(y,r)$ in positions $(1,2)$ and $(1,m+1)$.
\end{enumerate}
These operations leave all previously cleared entries in columns $3,\ldots,m$ unchanged. At each iteration, the absolute value of the entry in position $(1,m+1)$ decreases from $|y|$ to $r<|y|$. The process therefore terminates with the entry in position $(1,m+1)$ equal to zero. The absolute value of the resulting entry in position $(1,2)$ is
\[
\gcd(x,y)=\gcd(a_{12},\ldots,a_{1m},a_{1(m+1)}),
\]
where the second equality uses the induction hypothesis and the fact that column $m+1$ was unchanged before this stage. Thus $P(m+1)$ holds, completing the induction.

By induction, $P(n)$ holds. Hence the first row of $A\Psi(w_n)$ is $(a_{11},a_{12}^{(n)},0,\ldots,0)$, with $|a_{12}^{(n)}|=L$. Therefore $a_{12}^{(n)}=\pm L$.
\begin{enumerate}[(i)]
\item If $a_{12}^{(n)}=L$, set $w=w_n$.
\item If $a_{12}^{(n)}=-L$, right-multiply by $\Psi(\iota_2)$ to change the sign of column $2$, and set $w=w_n\iota_2$.
\end{enumerate}
The resulting word $w \in \langle S'\rangle$ satisfies $g = \Psi(w) \in \langle \Psi(S') \rangle$, preserves column $1$ of $A$, and transforms the subrow vector $\mathbf{a}$ into $(L, 0, \dots, 0)$.

\end{proof}

We now prove Theorem~\ref{thm:gen_img_stab_e1}.
\begin{proof}
Let $M = (a_{ij}) \in \Psi(\Out_{F_{n_k}}(F_n,e_1^*))$. Our goal is to write $M$ as a product of elements of the generating set
\[\mathcal{S}^{(n)} = \kappa_n\big(\Gamma(k)\big) \,\cup\, \Psi(S)\]
and their inverses. We proceed by successively multiplying $M$ by suitable elements of $\langle\mathcal S^{(n)}\rangle$ so that the resulting product reduces to the identity matrix $I_{n}$.
We begin by examining the first row of $M$, namely
$M(1) = (a_{11}, a_{12}, \dots, a_{1n})$.

By Lemma~\ref{lem:subrow-reduction}, $(a_{12},\dots,a_{1n})$ can be reduced to $(L,0,\dots,0)$ using only $\{\mu_{ij} \mid 2 \le i, j \le n,\ i \neq j\} \cup \{\sigma_{ij} \mid 2\leq i<j\leq n\} \cup \{\iota_{2}\}\subset S$, where $L=kg$ and $g=\gcd(c_2,\ldots,c_n)$, since $a_{1i}=kc_i$ for some $c_i\in\mathbb{Z}$ and all $i\geq 2$. This means there exists an $E^{(1)}\in\langle\Psi(S')\rangle\subset\langle\mathcal{S}^{(n)}\rangle$ such that the first row of $M_1=ME^{(1)}$ is $M_1(1)=M(1)E^{(1)}=(a_{11},L,0,\dots,0)$, with $a_{11}$ completely unchanged.

Since $M\in\mathrm{GL}(n,\mathbb Z)$, its first row is primitive, and hence
\[
\gcd(a_{11},kg)=\gcd(a_{11},a_{12},\ldots,a_{1n})=1.
\]
Also, $a_{11}\equiv1\pmod{k}$ implies $\gcd(a_{11},k)=1$, so $\gcd(a_{11},k^2g)=1$. By B\'ezout's identity, there exist $\alpha_1,\alpha_2\in\mathbb Z$ such that
\[
\alpha_1a_{11}+\alpha_2k^2g=1.
\]

Now, take the embedding $\kappa_n : \mathrm{SL}(2, \mathbb{Z}) \hookrightarrow \mathrm{GL}(n, \mathbb{Z})$ and set $E^{(2)}=\kappa_n(E_2)$, where $E_2=\begin{pmatrix}
 \alpha_1 & - kg \\ \alpha_2k & a_{11}
\end{pmatrix}$. Note that $\det E_2=\alpha_1a_{11}+\alpha_2k^2g=1$, and since $\alpha_1a_{11}\equiv 1 \pmod{k}$ and $a_{11}\equiv 1 \pmod{k}$, we have $\alpha_1\equiv 1\pmod{k}$. Hence $E_2\in\Gamma(k)$ and $E^{(2)}\in\kappa_n\big(\Gamma(k)\big)$. We consider the element
\[E^{(2)} = \begin{pmatrix} \alpha_1 & - kg & 0 & \cdots & 0 \\ \alpha_2k & a_{11} & 0 & \cdots & 0 \\ 0 & 0 & 1 & \cdots & 0 \\ \vdots & \vdots & \vdots & \ddots & \vdots \\ 0 & 0 & 0 & \cdots & 1 \end{pmatrix}.\]

Let $M_2=ME^{(1)}E^{(2)}$. Since $E^{(2)}$ is block diagonal, the first-row calculation reduces to
\[
\begin{pmatrix}a_{11}&kg\end{pmatrix}
\begin{pmatrix}\alpha_1&-kg\\\alpha_2k&a_{11}\end{pmatrix}
=\begin{pmatrix}1&0\end{pmatrix}.
\]
Thus $M_2(1)=(1,0,\ldots,0)$.
By construction, the $(1,1)$-entry of $(ME^{(1)})E^{(2)}$, denoted by $M_2(1,1)$, is equal to $1$, and the remaining entries of the first row of $M_2$ are $0$. Now, let $b_2\in\mathbb{Z}$ be the $(2, 1)$-entry of $M_2$. We perform the row operation $R_2 \longrightarrow R_2 - b_2 R_1$, where $R_2$ and $R_1$ denote the second and first rows, respectively. Since the first row of $M_2$ is $(1,0,\dots,0)$, this makes the $(2,1)$-entry $0$ and leaves the other entries of $R_2$ unchanged. The corresponding matrix operation is left multiplication by $E^{(2,1)}=\Psi(\mu_{12}^{-b_2})$, and the resulting matrix is $M_3=E^{(2,1)}M_2$. Similarly, for $3\leq i\leq n$, let $b_i\in\mathbb{Z}$ be the $(i,1)$-entry of $M_2$ and let $R_i$ denote the $i$th row. We perform the row operations $R_i \longrightarrow R_i - b_i R_1$, which make the $(i,1)$-entry $0$, and the corresponding operations are left multiplication by $E^{(i,1)}=\Psi(\mu_{1i}^{-b_i})$.
Using these matrices, we now modify $M_2$ to $\widetilde M \in \Psi(\Out_{F_{n_k}}(F_n,e_1^*))$, whose first column is $(1,0,\dots,0)^T$.

Each $\mu_{1i}$ is conjugate to $\mu_{12}$ by a permutation of $x_2,\ldots,x_n$, namely $\mu_{1i}=\sigma_{2i}\mu_{12}\sigma_{2i}$ for $3\leq i\leq n$, so these elementary matrices belong to $\langle\Psi(S)\rangle$. Now, we have $E^{(1)},E^{(2)},E^{(2,1)},\dots,E^{(n,1)}\in\langle\mathcal{S}^{(n)}\rangle$ such that
$E^{(n,1)}\cdots E^{(2,1)}\,ME^{(1)}E^{(2)}=\widetilde M$. Since the first row of $\widetilde M$ is $(1,0,\ldots,0)$ and its first column is $(1,0,\ldots,0)^T$, we have $\widetilde{M}=\begin{pmatrix}1&\mathbf{0}^T\\\mathbf{0}&B\end{pmatrix}$, where $B\in\mathrm{GL}(n-1,\mathbb Z)$ because $\det B=\det\widetilde M=\pm 1$. Note that $\{\Psi(\mu_{ij})\mid 2\le i,j \le n,\ i\ne j\}\subseteq\Psi(S')$. Using these matrices together with $D=\Psi(\iota_2)$ and all $\Psi(\sigma_{ij})$, ordinary integer Gaussian elimination reduces $B$ to $I_{n-1}$, and hence $\widetilde{M}$ to $I_n$, in finitely many operations.

Hence $\widetilde M\in\langle\Psi(S')\rangle$. Reversing the preceding operations gives
\[
M=\bigl(E^{(2,1)}\bigr)^{-1}\cdots\bigl(E^{(n,1)}\bigr)^{-1}
\widetilde M\bigl(E^{(2)}\bigr)^{-1}\bigl(E^{(1)}\bigr)^{-1}
\in\langle\mathcal S^{(n)}\rangle,
\]
which proves the assertion.
\end{proof}

\subsection{Generating the kernel of the restricted homology representation}
We begin by recalling the generating set discussed earlier in \eqref{eq:stabilizer-image-generators} and Corollary~\ref{cor:coset-representatives}:
\[S=\{\mu_{ij} \mid 2 \le i,j \le n,\ i \neq j\}\cup \{\sigma_{ij} \mid 2 \le i<j \le n\}\cup \{\iota_2,\mu_{12}\}.\]
Moreover,
$\Out_{F_{n_k}}(F_n)/\Out_{F_{n_k}}(F_n,e_1^*)$
has a set of coset representatives given by
\[\mathcal{S}''_{k}=\left\{\mu_{12}^{1-\bar{\ell}}\mu_{21}^{-1}\mu_{12}^{1-\ell}\mu_{21}^{\bar{\ell}}:\ell\in\mathbb{Z}_k^\times,\;\ell\bar{\ell}\equiv 1\pmod{k}
\right\}.\]
Recall the set $\widetilde{\mathcal S}_k$ from Section~\ref{sec:intro}. We next recall a theorem of Magnus~\cite{Magnus1935}, which gives a finite generating set for the Torelli subgroup $\mathcal{T}_n$.

\begin{theorem}[Magnus~\cite{Magnus1935}]
Let $F_n=\langle x_1,\ldots,x_n\rangle$. The group $\mathcal{T}_n$ is generated by the $\Out(F_n)$-classes of the automorphisms
\[\left\{N_{ij}:x_i\mapsto x_jx_ix_j^{-1}\mid i\neq j
\right\}\cup\left\{N_{ij\ell}:x_i\mapsto x_i[x_j,x_\ell]\mid i,j,\ell\text{ distinct},\ j<\ell\right\}.\]
Here, the automorphisms fix $x_m$ for $m\neq i$, and $[x_j,x_\ell]:=x_jx_\ell x_j^{-1}x_\ell^{-1}$ denotes the commutator of $x_j$ and $x_\ell$.
\end{theorem}

Recall that
\[
\mathcal N_n=\{N_{1j},N_{j1}\mid2\leq j\leq n\}
\cup\{N_{1j\ell}\mid2\leq j<\ell\leq n\}
\cup\{N_{i\ell}\mid2\leq i,\ell\leq n,\ i\neq\ell\}.
\]
In particular, $\mathcal N_n$ contains every basis-conjugating automorphism $N_{ij}$. To prove the required kernel containment, we express the remaining Magnus generators using these automorphisms and the Nielsen transformations already present in the proposed generating set.
\begin{theorem}\label{thm:kernel_Psi}
Let $k \ge 2$ and $n \ge 3$. Then
\[\ker\!\big(\Psi|_{\Out_{F_{n_k}}(F_n,e_1^*)}\big)\subset\left\langle
\mathcal{N}_n\cup\{\iota_2,\mu_{23}\}\cup\{\mu_{1i}\mid 2\leq i\leq n\}
\cup\{\sigma_{ij}\mid 2\leq i<j\leq n\}\right\rangle.\]
\end{theorem}
\begin{proof}
Let $G_n$ denote the subgroup on the right-hand side of the asserted containment. Every element of $\mathcal T_n=\ker\Psi$ fixes $e_1^*$ modulo $k$, so
$\ker\bigl(\Psi|_{\Out_{F_{n_k}}(F_n,e_1^*)}\bigr)=\mathcal T_n$.
It therefore suffices to show that every Magnus generator belongs to $G_n$.

All the generators $N_{ij}$ belong to $\mathcal N_n\subseteq G_n$. We next observe that $\mu_{i\ell}\in G_n$ whenever $\ell\geq2$ and $i\neq\ell$. For $i=1$, this is part of the definition of $G_n$. For $i\geq2$, choose a permutation $\pi$ of $\{2,\ldots,n\}$ with $\pi(2)=i$ and $\pi(3)=\ell$, extended to fix $1$. Such a permutation exists because $n\geq3$ and $i\neq\ell$. It is generated by the transpositions listed in $G_n$, and $\pi\mu_{23}\pi^{-1}=\mu_{i\ell}$.

Here and below, products of automorphisms mean composition, with the rightmost map applied first. For any three distinct indices $i,j,\ell$, the following identity holds in $\Aut(F_n)$:
\begin{equation}\label{eq:magnus-nielsen-identity}
N_{ij\ell}=N_{ij}^{-1}\mu_{i\ell}N_{ij}\mu_{i\ell}^{-1}.
\end{equation}
Indeed, both sides fix $x_m$ for $m\neq i$, while
\[
\begin{aligned}
\bigl(N_{ij}^{-1}\mu_{i\ell}N_{ij}\mu_{i\ell}^{-1}\bigr)(x_i)
&=N_{ij}^{-1}\bigl(x_jx_ix_\ell x_j^{-1}x_\ell^{-1}\bigr)\\
&=x_ix_jx_\ell x_j^{-1}x_\ell^{-1}
=x_i[x_j,x_\ell].
\end{aligned}
\]
Thus the same identity holds after passing to outer classes. In the Magnus generating set we have $j<\ell$, hence $\ell\geq2$. Both $N_{ij}$ and $\mu_{i\ell}$ therefore belong to $G_n$, so \eqref{eq:magnus-nielsen-identity} gives $N_{ij\ell}\in G_n$. Magnus's theorem now implies $\mathcal T_n\subseteq G_n$, as required.
\end{proof}

\begin{cor}\label{cor:314}
Let $k \ge 2$ and $n \ge 3$. Then
\[\Out_{F_{n_k}}(F_n,e_1^*)=\left\langle\widetilde{\mathcal{S}}_k\cup\mathcal{N}_n\cup\{\iota_2,\mu_{23}\}\cup\{\mu_{1i}\mid 2\leq i\leq n\}
\cup\{\sigma_{ij}\mid 2\leq i<j\leq n\}\right\rangle.\]
\end{cor}
\begin{proof}
Let $L_{n,k}$ denote the subgroup displayed on the right. Each of its generators fixes $e_1^*$ modulo $k$, so $L_{n,k}\leq\Out_{F_{n_k}}(F_n,e_1^*)$. Theorem~\ref{thm:kernel_Psi} gives $\mathcal T_n\subseteq L_{n,k}$. Moreover, $L_{n,k}$ contains the generators in $S$: the lower-block Nielsen transformations are conjugates of $\mu_{23}$ by lower-block permutations, and $\mu_{12}$ is among the listed $\mu_{1i}$. Together with $\widetilde{\mathcal S}_k$, their homology images generate the full image of the stabilizer by Theorem~\ref{thm:gen_img_stab_e1}. For any $\phi$ in the stabilizer, choose $g\in L_{n,k}$ with $\Psi(g)=\Psi(\phi)$. Then $\phi g^{-1}\in\mathcal T_n\subseteq L_{n,k}$, so $\phi\in L_{n,k}$.
\end{proof}

Combining Corollary~\ref{cor:coset-representatives} with  Corollary~\ref{cor:314}, we arrive at our main result.
\begin{theorem}[Main theorem]\label{thm:main-generating-set}
Let $k \ge 2$ and $n \ge 3$. Then
\[\Out_{F_{n_k}}(F_n)=\left\langle \mathcal S_k''\cup \widetilde{\mathcal{S}}_k\cup\mathcal{N}_n\cup\{\iota_2,\mu_{23}\}\cup\{\mu_{1i}\mid 2\leq i\leq n\}\cup\{\sigma_{ij}\mid 2\leq i<j\leq n\}\right\rangle.\]
\end{theorem}
\section{Lifting outer automorphism classes to the infinite chain of roses}\label{sec:Inf_cover}

Recall the infinite chain of roses $\mathcal C$ from Section~\ref{sec:intro}. Define the unit translation $\tau :\mathcal{C}\rightarrow\mathcal{C}$ by shifting the index of all cells by -1, which is $\tau(v_n)=v_{n-1}$, $\tau(e_{1,n})=e_{1,n-1}$, and $\tau(a_i^n)=a_i^{n-1}$, where $\langle \tau \rangle \cong \mathbb{Z}$.  In this section, we examine the liftable outer automorphism classes of the regular infinite-sheeted cyclic cover
$q_n:\mathcal{C}\to R_n$ induced by a $\mathbb Z$-action on $\mathcal{C}$, where the $\mathbb Z$-action is generated by the unit translation on $\mathcal C$ (see Figure~\ref{fig:inf-cover}). As in the finite-sheeted case, $(q_n)_\# H_1(\mathcal{C})\triangleleft H_1(R_n)$ and is identified with the kernel of the map $\bar \varphi:H_1(R_n)\to\mathbb{Z}$ given by $v\mapsto \langle a_1,v\rangle$. Consider the induced homomorphism $(q_n)_\#:H_1(\mathcal{C},\mathbb{Z})\to H_1(R_n,\mathbb{Z})$ given by
\begin{equation}
\begin{aligned}
a_i^k &\xmapsto{(q_n)_\#} \bar{a}_i,
&&\text{for } k\in\mathbb{Z},\ 2\leq i\leq n.
\end{aligned}\label{eq:covering-map}
\end{equation}
\noindent We use the techniques of Dowdall and Taylor \cite{spence} to show that the cover $q_n$ has an analogous Birman--Hilden short exact sequence.
\begin{lemma}\label{lem:infinite-cover-sequence}
Let $F_n=\langle x_1,\ldots,x_n\rangle$ with $n\geq2$, and let $F_\infty$ be the normal subgroup corresponding to the infinite-sheeted cyclic cover $q_n:\mathcal C\to R_n$. Thus $F_\infty$ is the kernel of the homomorphism $\varphi:F_n\to\mathbb Z$ given by $x_1\mapsto1$ and $x_i\mapsto0$ for $i\geq2$, and it is free of countably infinite rank. There is a short exact sequence
\begin{equation}\label{eq:seq_inf}
1 \longrightarrow \mathbb{Z} \longrightarrow \Out^{F_n}(F_\infty)\longrightarrow \Out_{F_\infty}(F_n) \longrightarrow 1.
\end{equation}
\end{lemma}

\begin{proof}
We define $\Aut_{F_\infty}(F_n)$ to be the subgroup of $\Aut(F_n)$ consisting of automorphisms that preserve $F_\infty$. The restriction of these automorphisms to $F_\infty$ determines a subgroup $\Aut^{F_n}(F_\infty)$ of $\Aut(F_\infty)$, given by the image of the map $\operatorname{res}: \Aut_{F_\infty}(F_n) \rightarrow \Aut(F_\infty)$. We further define $\Out_{F_\infty}(F_n)$ and $\Out^{F_n}(F_\infty)$ as the images of $\Aut_{F_\infty}(F_n)$ and $\Aut^{F_n}(F_\infty)$ under the respective quotient maps $p_{F_n}: \Aut(F_n) \rightarrow \Out(F_n)$ and $p_{F_\infty}: \Aut(F_\infty) \rightarrow \Out(F_\infty)$.

The restriction map $\operatorname{res}: \Aut_{F_\infty}(F_n) \rightarrow \Aut^{F_n}(F_\infty)$ is surjective by definition of $\Aut^{F_n}(F_\infty)$. It is also injective. To see this, suppose $\phi_1, \phi_2 \in \Aut_{F_\infty}(F_n)$ have the same restriction to $F_\infty$. Then the automorphism $\phi = \phi_1^{-1}\phi_2 \in \Aut_{F_\infty}(F_n)$ restricts to the identity on $F_\infty$, which means it fixes each element of $F_\infty$. Since the generators $x_i$ for $i \ge 2$ map to $0$ in $\mathbb{Z}$, they belong to $F_\infty$. Hence, $\phi(x_i) = x_i$ for all $i \ge 2$. Furthermore, the conjugate $x_1 x_2 x_1^{-1}$ maps to $0$ in $\mathbb{Z}$, so $x_1 x_2 x_1^{-1} \in F_\infty$. Thus,
\[
\begin{aligned}
\phi(x_1x_2x_1^{-1})&=x_1x_2x_1^{-1},\\
\phi(x_1)x_2\phi(x_1)^{-1}&=x_1x_2x_1^{-1}.
\end{aligned}
\]
which gives $(x_1^{-1}\phi(x_1)) x_2 (x_1^{-1}\phi(x_1))^{-1} = x_2$, implying that $x_1^{-1}\phi(x_1)$ commutes with $x_2$. The centralizer of a nontrivial element of a free group that is not a proper power is the infinite cyclic subgroup generated by that element~\cite{LyndonSchupp,MagnusKarrassSolitar}.
Since $x_2$ is a free basis element, it is not a proper power, so
$x_1^{-1}\phi(x_1) = x_2^k$ for some $k \in \mathbb{Z}$, that is,
$\phi(x_1)=x_1x_2^k$. Applying $\phi$ to $x_1^2x_2x_1^{-2}\in F_\infty$ gives
\begin{align*}
\phi(x_1^2 x_2 x_1^{-2}) &= x_1^2 x_2 x_1^{-2},\\
(x_1 x_2^k)^2 x_2 (x_1 x_2^k)^{-2} &= x_1^2 x_2 x_1^{-2},\\
x_1 x_2^k x_1 x_2 x_1^{-1} x_2^{-k} x_1^{-1} &= x_1^2 x_2 x_1^{-2}.
\end{align*}
It follows that $x_2^k$ must commute with $x_1 x_2 x_1^{-1}$. The elements $x_2$ and $x_1x_2x_1^{-1}$ freely generate a rank-two subgroup, so $x_2^k$ cannot commute with $x_1x_2x_1^{-1}$ unless $k=0$. Therefore, $\phi(x_1) = x_1$, which means $\operatorname{res}$ is an isomorphism.
As $F_\infty$ is a normal subgroup of $F_n$, conjugation by any element of $F_n$ preserves $F_\infty$. This means $\operatorname{Inn}(F_n) \subseteq \Aut_{F_\infty}(F_n)$. Therefore, the kernel of the projection $p_{F_n}$ is exactly $\Aut_{F_\infty}(F_n) \cap \operatorname{Inn}(F_n) = \operatorname{Inn}(F_n)$. Thus, the kernel of the composition $\hat{\pi} := p_{F_n} \circ \operatorname{res}^{-1}: \Aut^{F_n}(F_\infty) \rightarrow \Out_{F_\infty}(F_n)$ is exactly $\operatorname{Inn}(F_n)$.

Observe that every inner automorphism of $F_\infty$ lifts to an inner automorphism of $F_n$ via the map $\operatorname{res}^{-1}$. Since the resulting inner automorphism of $F_n$ clearly normalizes $F_\infty$, we have $\ker(p_{F_\infty}) \subseteq\ker(\hat{\pi})$.
It follows that the homomorphism descends to a well-defined epimorphism $\pi: \Out^{F_n}(F_\infty) \rightarrow \Out_{F_\infty}(F_n)$ with kernel $\ker(\pi) = \ker(\hat{\pi})/\ker(p_{F_\infty}) \cong F_n/F_\infty$. Thus we recover the short exact sequence
\begin{equation*}
1 \longrightarrow \mathbb{Z} \longrightarrow \Out^{F_n}(F_\infty)\longrightarrow \Out_{F_\infty}(F_n) \longrightarrow 1.
\end{equation*}
\end{proof}
Define
\[
\operatorname{UOut}(F_n):=\bigcap_{k\geq2}\Out_{F_{n_k}}(F_n).
\]
The Torelli subgroup $\mathcal T_n$ lies in $\operatorname{UOut}(F_n)$. Theorem~\ref{thm:homo_criteria} gives the following integral criterion.
\begin{cor}\label{cor:UOut}
Let $\phi \in \Out(F_n)$ and let $\Psi(\phi) = (c_{ij})_{n \times n}$.
Then $\phi \in \operatorname{UOut}(F_n)$ if and only if $c_{1i}=0$ for $2\leq i\leq n$ and $c_{11}=\pm1$.
\end{cor}
\noindent Using this criterion, we establish the following.
\begin{prop}\label{prop:uout_F_inf}
For $n \geq 2$, $\operatorname{UOut}(F_n)=\Out_{F_\infty}(F_n)$.
\end{prop}
\begin{proof}
From the definition of $(q_n)_\#$ in Equation~\ref{eq:covering-map}, we have
\[W := (q_n)_\#(H_1(\mathcal{C},\mathbb{Z}))
= \langle \bar a_2,\bar a_3,\ldots,\bar a_{n}\rangle.
\]
An outer automorphism lifts to $\mathcal C$ precisely when its induced automorphism of $H_1(R_n,\mathbb Z)$ preserves $W$, since the covering subgroup is the kernel of the homomorphism $F_n\to\mathbb Z$ that records the exponent sum of $x_1$. The subgroup $W$ is the kernel of the integral covector $e_1^*$. An invertible integral matrix preserves $W$ if and only if its first row is $(\pm1,0,\ldots,0)$. By Corollary~\ref{cor:UOut}, this is exactly the condition for membership in $\operatorname{UOut}(F_n)$, proving the assertion.
\end{proof}
It is easy to observe that
$\bigcap\limits_{k\geq 2} \Out_{F_{n_k}}(F_n,e_1^*)
= \Stab_{\Out(F_n)}(e_1^*)$
is an index $2$ subgroup of $\operatorname{UOut}(F_n)$, where the $e_1^*$
appearing in $\Stab_{\Out(F_n)}(e_1^*)$ refers to
$e_1^*\in H^1(R_n,\mathbb{Z})$. Combining this index-two description with Theorem~\ref{thm:kernel_Psi}, we obtain the following corollary.
\begin{cor}\label{cor:gen_f_inf}
For $n\geq 3$,
\[\operatorname{UOut}(F_n)=\left\langle\mathcal{N}_n\cup\{\iota_2,\mu_{23},\iota\}\cup\{\mu_{1i}\mid 2\leq i\leq n\}\cup\{\sigma_{ij}\mid 2\leq i<j\leq n\}\right\rangle.
\]
\end{cor}
\begin{proof}
Let $G_n$ be the subgroup on the right-hand side of Theorem~\ref{thm:kernel_Psi}, and put $K=\Stab_{\Out(F_n)}(e_1^*)$, where $e_1^*\in H^1(R_n,\mathbb Z)$. Every generator of $G_n$ fixes this integral covector, so $G_n\leq K$. The matrices $\Psi(\mu_{1i})=I+E_{i1}$ generate all integral lower-left columns. The images of $\mu_{23}$, $\iota_2$, and the lower-block permutations generate the entire lower-right $\mathrm{GL}(n-1,\mathbb Z)$ block. Consequently,
\[
\Psi(G_n)=\Psi(K)=
\left\{\begin{pmatrix}1&0\\v&B\end{pmatrix}:
 v\in\mathbb Z^{n-1},\ B\in\mathrm{GL}(n-1,\mathbb Z)\right\}.
\]
By Theorem~\ref{thm:kernel_Psi}, $\mathcal T_n\subseteq G_n$. Equality of the homological images therefore implies $G_n=K$. Finally, $K$ has index $2$ in $\operatorname{UOut}(F_n)$, and $\iota\in\operatorname{UOut}(F_n)\setminus K$. Hence $\operatorname{UOut}(F_n)=\langle G_n,\iota\rangle$, proving the asserted generating statement.
\end{proof}

In Section~\ref{sec:intro}, we stated that $\Out_{F_\infty}(F_n)$ is the inverse limit of the liftable subgroups $\Out_{F_{n_k}}(F_n)$ of the finite-sheeted covers $p_k$. To make the limit description precise, let $D=\{k\in\mathbb{N}:k\geq 2\}$, ordered by divisibility, and put
\[
H_k=\Out_{F_{n_k}}(F_n),\qquad K_k=\Out_{F_{n_k}}(F_n,e_1^*).
\]
If $k\mid\ell$, the homological criteria give inclusions $H_\ell\hookrightarrow H_k$ and $K_\ell\hookrightarrow K_k$. With respect to these inclusions, $\{H_k\}_{k\in D}$ and $\{K_k\}_{k\in D}$ form two inverse systems over $D$. An element of either inverse limit is a compatible family $(g_k)_{k\in D}$, that is, $g_\ell=g_k$ whenever $k\mid\ell$. Such a family has equal entries at any two indices $k$ and $\ell$, since both entries agree with the entry at $\operatorname{lcm}(k,\ell)$. Consequently, the inverse limits identify canonically with the intersections inside $\Out(F_n)$:
\[
\begin{aligned}
\varprojlim_{k\in D} H_k &= \bigcap_{k\geq 2} H_k = \operatorname{UOut}(F_n),\\
\varprojlim_{k\in D} K_k &= \bigcap_{k\geq 2} K_k = \Stab_{\Out(F_n)}(e_1^*)=\Out(F_n,e_1^*),
\end{aligned}
\]
where $e_1^*\in H^1(R_n,\mathbb{Z})$ in the last stabilizer. Since $\Out(F_n,e_1^*)$ has index $2$ in $\operatorname{UOut}(F_n)$ and $\iota\in\operatorname{UOut}(F_n)\setminus\Out(F_n,e_1^*)$, we obtain the following normal series.
\begin{cor}
For $n\geq3$, the group $\operatorname{UOut}(F_n)$ admits the normal series
\[1\triangleleft\mathcal{T}_n\triangleleft\Out(F_n,e_1^*)
\triangleleft\operatorname{UOut}(F_n),
\]
where the two distinct cosets of
$\operatorname{UOut}(F_n)/\Out(F_n,e_1^*)$ are represented by
the elements in $\langle\iota\rangle$.
\end{cor}
\section{Acknowledgments}
The first author is supported by the Prime Minister's Research Fellowship (ID: 0403022), Ministry of Education, Government of India. The third author is supported by an ANRF MATRICS grant (File No. MTR/2023/000610), Department of Science and Technology, India.

\section{Tool use} In preparing this work, the authors used ChatGPT for language editing and grammar refinement. After utilizing this tool, the authors thoroughly reviewed and edited the content as necessary and assume full responsibility for its content.

\bibliographystyle{plain}
\bibliography{Liftable_out} I
\end{document}